# Estimation in a class of semiparametric transformation models

**Dorota M. Dabrowska**[1],*

*University of California, Los Angeles*

**Abstract:** We consider estimation in a class of semiparametric transformation models for right–censored data. These models gained much attention in survival analysis; however, most authors consider only regression models derived from frailty distributions whose hazards are decreasing. This paper considers estimation in a more flexible class of models and proposes conditional rank M-estimators for estimation of the Euclidean component of the model.

## 1. Introduction

Semiparametric transformation models provide a common tool for regression analysis. We consider estimation in a class of such models designed for analysis of failure time data with time independent covariates. Let $\mu$ be the marginal distribution of a covariate vector $Z$ and let $H(t|z)$ be the cumulative hazard function of the conditional distribution of failure time $T$ given $Z$. We assume that for $\mu$–almost all $z$ ($\mu$ a.e. $z$) this function is of the form

$$(1.1) \qquad H(t|z) = A(\Gamma(t), \theta|z)$$

where $\Gamma$ is an unknown continuous increasing function mapping the support of the failure time $T$ onto the positive half-line. For $\mu$ a.e. $z$, $A(x, \theta|z)$ is a conditional cumulative hazard function dependent on a Euclidean parameter $\theta$ and having hazard rate $\alpha(x, \theta|z)$ strictly positive at $x = 0$ and supported on the whole positive half-line. Special cases include

(i) the proportional hazards model with constant hazard rate $\alpha(x, \theta|z) = \exp(\theta^T z)$ (Lehmann [23], Cox [12]);
(ii) transformations to distributions with monotone hazards such as the proportional odds and frailty models or linear hazard rate regression model (Bennett [2], Nielsen *et al.* [28], Kosorok *et al.* [22], Bogdanovicius and Nikulin [9]);
(iii) scale regression models induced by half-symmetric distributions (section 3).

The proportional hazards model remains the most commonly used transformation model in survival analysis. Transformation to exponential distribution entails that for any two covariate levels $z_1$ and $z_2$, the ratio of hazards is constant in $x$ and equal to $\alpha(x, \theta|z_1)/\alpha(x, \theta|z_2) = \exp(\theta^T[z_1 - z_2])$. Invariance of the model with respect to monotone transformations enstails that this constancy of hazard ratios is preserved by the transformation model. However, in many practical circumstances

---

*Research supported in part by NSF grant DMS 9972525 and NCI grant 2R01 95 CA 65595-01.
[1]Department of Biostatistics, School of Public Health, University of California, Los Angeles, CA 90095-1772, e-mail: dorota@ucla.edu

this may fail to hold. For example, a new treatment ($z_1 = 1$) may be initially beneficial as compared to a standard treatment ($z_2 = 0$), but the effects may decay over time, $\alpha(x, \theta|z_1 = 1)/\alpha(x, \theta|z_2 = 0) \downarrow 1$ as $x \uparrow \infty$. In such cases the choice of the proportional odds model or a transformation model derived from frailty distributions may be more appropriate. On the other hand, transformation to distributions with increasing or non-monotone hazards allows for modeling treatment effects which have divergent long-term effects or crossing hazards. Transformation models have also found application in regression analyses of multivariate failure time data, where models are often defined by means of copula functions and marginals are specified using models (1.1).

We consider parameter estimation in the presence of right censoring. In the case of uncensored data, the model is invariant with respect to the group of increasing transformations mapping the positive half-line onto itself so that estimates of the parameter $\theta$ are often sought within the class of conditional rank statistics. Except for the proportional hazards model, the conditional rank likelihood does not have a simple tractable form and estimation of the parameter $\theta$ requires joint estimation of the pair $(\theta, \Gamma)$. An extensive study of this estimation problem was given by Bickel [4], Klaassen [21] and Bickel and Ritov [5]. In particular, Bickel [4] considered the two sample testing problem, $\mathcal{H}_0 : \theta = \theta_0$ vs $\mathcal{H} : \theta > \theta_0$, in one-parameter transformation models. He used projection methods to show that a nonlinear rank statistic provides an efficient test, and applied Sturm-Liouville theory to obtain the form of its score function. Bickel and Ritov [5] and Klaassen [21] extended this result to show that under regularity conditions, the rank likelihood in regression transformation models forms a locally asymptoticaly normal family and estimation of the parameter $\theta$ can be based on a one-step MLE procedure, once a preliminary $\sqrt{n}$ consistent estimate of $\theta$ is given. Examples of such estimators, specialized to linear transformation models, can be found in [6, 13, 15], among others.

In the case of censored data, the estimation problem is not as well understood. Because of the popularity of the proportional hazards model, the most commonly studied choice of (1.1) corresponds to transformation models derived from frailty distributions. Murphy *et al.* [27] and Scharfstein *et al.* [31] proposed a profile likelihood method of analysis for the generalized proportional odds ratio models. The approach taken was similar to the classical proportional hazards model. The model (1.1) was extended to include all monotone functions $\Gamma$. With fixed parameter $\theta$, an approximate likelihood function for the pair $(\theta, \Gamma)$ was maximized with respect to $\Gamma$ to obtain an estimate $\Gamma_{n\theta}$ of the unknown transformation. The estimate $\Gamma_{n\theta}$ was shown to be a step function placing mass at each uncensored observation, and the parameter $\theta$ was estimated by maximizing the resulting profile likelihood. Under certain regularity conditions on the censoring distribution, the authors showed that the estimates are consistent, asymptotically Gaussian at rate $\sqrt{n}$, and asymptotically efficient for estimation of both components of the model. The profile likelihood method discussed in these papers originates from the counting process proportional hazards frailty intensity models of Nielsen et al. [28]. Murphy [26] and Parner [30] developed properties of the profile likelihood method in multi-jump counting process models. Kosorok *et al* [22] extended the results to one-jump frailty intensity models with time dependent covariates, including the gamma, the lognormal and the generalized inverse Gaussian frailty intensity models. Slud and Vonta [33] provided a separate study of consistency properties of the nonparametric maximum profile likelihood estimator in transformation models assuming that the cumulative hazard function (1.1) is of the form $H(t|z) = A(\exp[\theta^T z]\Gamma(t))$ where $A$ is a known concave function.

Several authors proposed also ad hoc estimates of good practical performance. In particular, Cheng *et al.* [11] considered estimation in the linear transformation model in the presence of censoring independent of covariates. They showed that estimation of the parameter $\theta$ can be accomplished without estimation of the transformation function by means of U-statistics estimating equations. The approach requires estimation of the unknown censoring distribution, and does not extend easily to models with censoring dependent on covariates. Further, Yang and Prentice [34] proposed minimum distance estimation in the proportional odds ratio model and showed that the unknown odds ratio function can be estimated based on a sample analogue of a linear Volterra equation. Bogdanovicius *et al.* [9, 10] considered estimation in a class of generalized proportional hazards intensity models that includes the transformation model (1.1) as a special case and proposed a modified partial likelihood for estimation of the parameter $\theta$. As opposed to the profile likelihood method, the unknown transformation was profiled out from the likelihood using a martingale-based estimate of the unknown transformation obtained by solving recurrently a Volterra equation.

In this paper we consider an extension of estimators studied by Cuzick [13] and Bogdanovicius *et al.* [9, 10] to a class of M-estimators of the parameter $\theta$. In Section 2 we shall apply a general method for construction of M-estimates in semiparametric models outlined in Chapter 7 of Bickel *et al.* [6]. In particular, the approach requires that the nuisance parameter and a consistent estimate of it be defined in a larger model $\mathcal{P}$ than the stipulated semiparametric model. Denoting by $(X, \delta, Z)$, the triple corresponding to a nonnegative time variable $X$, a binary indicator $\delta$ and a covariate $Z$, in this paper we take $\mathcal{P}$ as the class of all probability measures such that the covariate $Z$ is bounded and the marginal distribution of the withdrawal times is either continuous or has a finite number of atoms. Under some regularity conditions on the core model $\{A(x, \theta|z) : \theta \in \Theta, x > 0\}$, we define a parameter $\Gamma_{P,\theta}$ as a mapping of $\mathcal{P} \times \Theta$ into a convex set of monotone functions. The parameter represents a transformation function that is defined as a solution to a nonlinear Volterra equation. We show that its "plug-in" estimate $\Gamma_{P_n,\theta}$ is consistent and asymptotically linear at rate $\sqrt{n}$. Here $P_n$ is the empirical measure of the data corresponding to an iid sample of the $(X, \delta, Z)$ observations. Further, we propose a class of M-estimators for the parameter $\theta$. The estimate will be obtained by solving a score equation $U_n(\theta) = 0$ or $U_n(\theta) = o_P(n^{-1/2})$ for $\theta$. Similarly to the case of the estimator $\Gamma_{n\theta}$, the score function $U_n(\theta)$ is well defined (as a statistic) for any $P \in \mathcal{P}$. It forms, however, an approximate V-process so that its asymptotic properties cannot be determined unless the "true" distribution $P \in \mathcal{P}$ is defined in sufficient detail (Serfling [32]). The properties of the score process will be developed under the added assumption that at true $P \in \mathcal{P}$, the observation $(X, \delta, Z) \sim \mathcal{P}$ has the same distribution as $(T \wedge \tilde{T}, 1(T \leq \tilde{T}), Z)$, where $T$ and $\tilde{T}$ represent failure and censoring times conditionally independent given the covariate $Z$, and the conditional distribution of the failure time $T$ given $Z$ follows the transformation model (1.1).

Under some regularity conditions, we show that the M-estimates converge at rate $\sqrt{n}$ to a normal limit with a simple variance function. By solving a Fredholm equation of second kind, we also show that with an appropriate choice of the score process, the proposed class of censored data rank statistics includes estimators of the parameter $\theta$ whose asymptotic variance is equal to the inverse of the asymptotic variance of the M-estimating score function $\sqrt{n}U_n(\theta_0)$. We give a derivation of the resolvent and solution of the equation based on Fredholm determinant formula. We also show that this is a Sturm-Liouville equation, though of a different form than in [4, 5] and [21].

The class of transformation models considered in this paper is different than in the literature on nonparametric maximum likelihood estimation (NMPLE); in particular, hazard rates of core models need not be decreasing. In section 2, the core models are assumed to have hazards $\alpha(x, \theta, z)$ uniformly bounded between finite positive constants. With this aid we show that the mapping $\Gamma_{P,\theta}$ of $\mathcal{P} \times \Theta$ into the class of monotone functions is well defined on the entire support of the withdrawal time distribution, and without any special conditions on the probability distribution $P$. Under the assumption that the upper support point $\tau_0$ of the withdrawal time distribution is a discontinuity point, the function $\Gamma_{P,\theta}$ is shown to be bounded. If $\tau_0$ is a continuity point of this distribution, the function $\Gamma_{P,\theta}(t)$ is shown to grow to infinity as $t \uparrow \tau_0$. In the absence of censoring, the model (1.1) assumes that the unknown transformation is an unbounded function, so we require $\Gamma_{P,\theta}$ to have this property as well. In section 3, we use invariance properties of the model to show that the results can also be applied to hazards $\alpha(x, \theta, z)$ which are positive at the origin, but only locally bounded and locally bounded away from 0. All examples in this section refer to models whose conditional hazards are hyperbolic, i.e can be bounded (in a neighbourhood of the true parameter) between a linear function $a + bx$ and a hyperbola $(c + dx)^{-1}$, for some $a > 0, c > 0$ and $b \geq 0, d \geq 0$. As an example, we discuss the linear hazard rate transformation model, whose conditional hazard function is increasing, but its conditional density is decreasing or non-monotone, and the gamma frailty model with fixed frailty parameter or frailty parameters dependent on covariates.

We also examine in some detail scale regression models whose core models have cumulative hazards of the form $A_0(x \exp[\beta^T z])$. Here $A_0$ is a known cumulative hazard function of a half-symmetric distribution with density $\alpha_0$. Our results apply to such models if for some fixed $\xi \in [-1, 1]$ and $\eta \geq 0$, the ratio $\alpha_0/g, g(x) = [1 + \eta x]^\xi$ is a function locally bounded and locally bounded away from zero. We show that this choice includes half-logistic, half-normal and half-t scale regression models, whose conditional hazards are increasing or non-monotone while densities are decreasing. We also give examples of models (with coefficient $\xi \notin [-1, 1]$) to which the results derived here cannot be applied.

Finally, this paper considers only the gamma frailty model with the frailty parameter fixed or dependent on covariates. We show, however, that in the case that the known transformation is the identity map, the gamma frailty regression model (frailty parameter independent of covariates) is not regular in its entire parameter range. When the transformation is unknown, and the parameter set restricted to $\eta \geq 0$, we show that the frailty parameter controls the shape of the transformation. We do not know at the present time, if there exists a class of conditional rank statistics which allows to estimate the parameter $\eta$, without any additional regularity conditions on the unknown transformation.

In Section 4 we summarize the findings of this paper and outline some open problems. The proofs are given in the remaining 5 sections.

## 2. Main results

We shall first give regularity conditions on the model (Section 2.1). The asymptotic properties of the estimate of the unknown transformation are discussed in Section 2.2. Section 2.3 introduces some additional notation. Section 2.4 considers estimation of the Euclidean component of the model and gives examples of M-estimators of this parameter.

### 2.1. The model

Throughout the paper we assume that $(X, \delta, Z)$ is defined on a complete probability space $(\Omega, \mathcal{F}, P)$, and represents a nonnegative withdrawal time $(X)$, a binary indicator $(\delta)$ and a vector of covariates $(Z)$. Set $N(t) = 1(X \leq t, \delta = 1)$, $Y(t) = 1(X \geq t)$ and let $\tau_0 = \tau_0(P) = \sup\{t : E_P Y(t) > 0\}$. We shall make the following assumption about the "true" probability distribution $P$.

**Condition 2.0.** $P \in \mathcal{P}$ where $\mathcal{P}$ is the class of all probability distributions such that

(i) The covariate $Z$ has a nondegenerate marginal distribution $\mu$ and is bounded: $\mu(|Z| \leq C) = 1$ for some constant $C$.
(ii) The function $E_P Y(t)$ has at most a finite number of discontinuity points, and $E_P N(t)$ is either continuous or discrete.
(iii) The point $\tau > 0$ satisfies $\inf\{t : E_P[N(t)|Z = z] > 0\} < \tau$ for $\mu$ a.e. $z$. In addition, $\tau = \tau_0$, if $\tau_0$ is an discontinuity point of $E_P Y(t)$, and $\tau < \tau_0$, if $\tau_0$ is a continuity point of $E_P Y(t)$.

For given $\tau$ satisfying Condition 2.0(iii), we denote by $\|\cdot\|_\infty$ the supremum norm in $\ell^\infty([0, \tau])$. The second set of conditions refers to the core model $\{A(\cdot, \theta|z) : \theta \in \Theta\}$.

**Condition 2.1.** (i) The parameter set $\Theta \subset R^d$ is open, and $\theta$ is identifiable in the core model: $\theta \neq \theta'$ iff $A(\cdot, \theta|z) \not\equiv A(\cdot, \theta'|z)$ $\mu$ a.e. $z$.
(ii) For $\mu$ almost all $z$, the function $A(\cdot, \theta|z)$ has a hazard rate $\alpha(\cdot, \theta|z)$. There exist constants $0 < m_1 < m_2 < \infty$ such that $m_1 \leq \alpha(x, \theta|z) \leq m_2$ for $\mu$ a.e. $z$ and all $\theta \in \Theta$.
(iii) The function $\ell(x, \theta, z) = \log \alpha(x, \theta, z)$ is twice continuously differentiable with respect to both $x$ and $\theta$. The derivatives with respect to $x$ (denoted by primes) and with respect to $\theta$ (denoted by dots) satisfy

$$|\ell'(x, \theta, z)| \leq \psi(x), \quad |\ell''(x, \theta, z)| \leq \psi(x),$$
$$|\dot{\ell}(x, \theta, z)| \leq \psi_1(x), \quad |\ddot{\ell}(x, \theta, z)| \leq \psi_2(x),$$
$$|g(x, \theta, z) - g(x', \theta', z)| \leq \max(\psi_3(x), \psi_3(x'))[|x - x'| + |\theta - \theta'|],$$

where $g = \ddot{\ell}, \dot{\ell}'$ and $\ell''$. Here $\psi$ is a constant or a continuous bounded decreasing function. The functions $\dot{\ell}, \ddot{\ell}$ and $\dot{\ell}'$ are locally bounded and $\psi_p, p = 1, 2, 3$ are continuous, bounded or strictly increasing and such that $\psi_p(0) < \infty$,

$$\int_0^\infty e^{-x} \psi_1^2(x) dx < \infty, \quad \int_0^\infty e^{-x} \psi_p(x) dx < \infty, \quad p = 2, 3.$$

To evaluate the score process for estimation of the parameter $\theta$, we shall use the following added assumption.

**Condition 2.2.** The true distribution $P$, $P \in \mathcal{P}$, is the same as that of $(X, \delta, Z) \sim (T \wedge \tilde{T}, 1(T \leq \tilde{T}), Z)$, where $T$ and $\tilde{T}$ represent failure and censoring times. The variables $T$ and $\tilde{T}$ are conditionally independent given $Z$. In addition

(i) The conditional cumulative hazard function of $T$ given $Z$ is of the form $H(t|z) = A(\Gamma_0(t), \theta_0|z)$ $\mu$ a.e. $z$, where $\Gamma_0$ is a continuous increasing function, and $A(x, \theta_0|z) = \int_0^x \alpha(u, \theta_0|z) du, \theta_0 \in \Theta$, is a cumulative hazard function with hazard rate $\alpha(u, \theta_0|z)$ satisfying Conditions 2.1.

(ii) If $\tau_0$ is a discontinuity point of the survival function $E_P Y(t)$, then $\tau_0 = \sup\{t : P(\tilde{T} \geq t) > 0\} < \sup\{t : P(T \geq t) > 0\}$. If $\tau_0$ is a continuity point of this survival function, then $\tau_0 = \sup\{t : P(T \geq t) > 0\} \leq \sup\{t : P(\tilde{T} \geq t) > 0\}$.

For $P \in \mathcal{P}$, let $A(t) = A_P(t)$ be given by

$$A(t) = \int_0^t \frac{EN_P(du)}{E_P Y(u)}. \tag{2.1}$$

If the censoring time $\tilde{T}$ is independent of covariates, then $A(t)$ reduces to the marginal cumulative hazard function of the failure time $T$, restricted to the interval $[0, \tau_0]$. Under Assumption 2.2 this parameter forms in general a function of the marginal distribution of covariates, and conditional distributions of both failure and censoring times. Nevertheless, we shall find it, and the associated Aalen–Nelson estimator, quite useful in the sequel. In particular, under Assumption 2.2, the conditional cumulative hazard function $H(t|z)$ of $T$ given $Z$ is uniformly dominated by $A(t)$. We have

$$A(t) = \int_0^t E[\alpha(\Gamma_0(u-), \theta_0, Z)|X \geq u]\Gamma_0(du)$$

and

$$\frac{H(dt|z)}{A(dt)} = \frac{\alpha(\Gamma_0(t-), \theta_0, z)}{E\alpha(\Gamma_0(t-), \theta_0, Z)|X \geq t)},$$

for $t \leq \tau(z) = \sup\{t : EY(t)|Z = z > 0\}$ and $\mu$ a.e. $z$. These identities suggest to define a parameter $\Gamma_{P,\theta}$ as solution to the nonlinear Volterra equation

$$\begin{aligned}
\Gamma_{P,\theta}(t) &= \int_0^t \frac{E_P N(du)}{E_P Y(u)\alpha(\Gamma_\theta(u-), \theta, Z)} \\
&= \int_0^t \frac{A_P(du)}{E_P \alpha(\Gamma_\theta(u-), \theta, Z)|X \geq u)},
\end{aligned} \tag{2.2}$$

with boundary condition $\Gamma_{P,\theta}(0-) = 0$. Because Conditions 2.2 are not needed to solve this equation, we shall view $\Gamma$ as a map of the set $\mathcal{P} \times \Theta$ into $\mathcal{X} = \cup\{\mathcal{X}(P) : P \in \mathcal{P}\}$, where

$$\mathcal{X}(P) = \{g : g \text{ increasing}, e^{-g} \in D(\mathcal{T}), g \ll E_P N, m_2^{-1} A_P \leq g \leq m_1^{-1} A_P\}$$

and $m_1, m_2$ are constants of Condition 2.1(iii). Here $D(\mathcal{T})$ denotes the space of right-continuous functions with left-hand limits, and we choose $\mathcal{T} = [0, \tau_0]$, if $\tau_0$ is a discontinuity point of the survival function $E_P Y(t)$, and $\mathcal{T} = [0, \tau_0)$, if it is a continuity point. The assumption $g \ll E_P N$ means that the functions $g$ in $\mathcal{X}(P)$ are absolutely continuous with respect to the sub-distribution function $E_P N(t)$. The monotonicity condition implies that they admit integral representation $g(t) = \int_0^t h(u) dE_P N(u)$ and $h \geq 0$, $E_P N$-almost everywhere.

### 2.2. Estimation of the transformation

Let $(N_i, Y_i, Z_i)$, $i = 1, \ldots, n$ be an iid sample of the $(N, Y, Z)$ processes. Set $S(x, \theta, t) = n^{-1} \sum_{i=1}^n Y_i(t)\alpha(x, \theta, Z_i)$ and denote by $\dot{S}, S'$ the derivatives of these processes with respect to $\theta$ (dots) and $x$ (primes) and let $s, \dot{s}, s'$ be the corresponding expectations. Suppressing dependence of the parameter $\Gamma_{P,\theta}$ on P, set

$$C_\theta(t) = \int_0^t \frac{EN(du)}{s^2(\Gamma_\theta(u-), \theta, u)}.$$

For $u \leq t$, define also

$$
\begin{aligned}
\mathcal{P}_\theta(u,t) &= \pi_{(u,t]}(1 - s'(\Gamma_\theta(w-), \theta, w)C_\theta(dw)), \\
&= \exp[-\int_u^t s'(\Gamma_\theta(w-), \theta, w)C_\theta(dw)] \quad \text{if } EN(t) \text{ is continuous,} \\
&= \prod_{u < w \leq t} [1 - s'(\Gamma_\theta(w-), \theta, w)C_\theta(dw)] \quad \text{if } EN(t) \text{ is discrete.}
\end{aligned}
$$
(2.3)

Finally, we follow Bogdanovicius and Nikulin [9], and use

$$\Gamma_{n\theta}(t) = \int_0^t \frac{N_\cdot(du)}{S(\Gamma_{n\theta}(u-), \theta, u)}, \quad \Gamma_{n\theta}(0-) = 0, \theta \in \Theta,$$

to estimate the unknown transformation. Here $N_\cdot = n^{-1}\Sigma_{i=1}^n N_i$.

**Proposition 2.1.** *Let $P \in \mathcal{P}$ be a distribution satisfying Conditions 2.0(i), (ii) and let $(X_i, \delta_i, Z_i)$, $i = 1, \ldots, n$ be an iid sample from this distribution. Suppose that Conditions 2.1 are fulfilled by the family $\{A(\cdot, \theta|z) : \theta \in \Theta\}$, and let $\tau$ be an arbitrary point such that Condition 2.0(iii) holds.*

(i) *Equation (2.2) has a unique locally bounded solution satisfying $0 < \Gamma_\theta(\tau_0) < \infty$ if $\tau_0 = \tau_0(P)$ is a discontinuity point of $EY_P(t)$ and $\lim_{t \uparrow \infty} \Gamma_\theta(t) = \infty$ if $\tau_0$ is a continuity point of this survival function. For any point $\tau$, the plug-in estimate $\{\Gamma_{n\theta}(t) : t \leq \tau, \theta \in \Theta\}$ satisfies $\sup_{\theta \in \Theta} \|\Gamma_{n\theta} - \Gamma_\theta\|_\infty \to 0$ a.s. In addition, if $\tau_0$ is a continuity point of $E_P Y(t)$, then $\sup\{|\exp(-\Gamma_\theta) - \exp(-\Gamma_{n\theta})|(t) : \theta \in \Theta, t \in \mathcal{T}\} = o_P(1)$.*
(ii) *The function $\Theta \ni \theta \to \{\Gamma_\theta(t) : t \in [0,\tau]\} \in \ell^\infty([0,\tau])$ is Fréchet differentiable with respect to $\theta$ and the derivative satisfies*

$$\dot{\Gamma}_\theta(t) = -\int_0^t \dot{s}(\Gamma_\theta(u-), \theta, u)C_\theta(du)\mathcal{P}_\theta(u,t).$$

*The estimate $\{\dot{\Gamma}_{n\theta}(t) : t \leq \tau, \theta \in \Theta\}$ satisfies $\sup_{\theta \in \Theta} \|\dot{\Gamma}_{n\theta} - \dot{\Gamma}_\theta\|_\infty \to 0$ a.s.*
(iii) *The process $\{\hat{W}(t,\theta) = \sqrt{n}[\Gamma_{n\theta} - \Gamma_\theta](t) : t \leq \tau, \theta \in \Theta\}$ converges weakly in $\ell^\infty([0,\tau] \times \Theta)$ to*

$$W(t,\theta) = R(t,\theta) - \int_{[0,t]} R(u-, \theta)C_\theta(du)V_\theta(u,t),$$

*where $V_\theta(u,t) = 1(u < t)s'(\Gamma_\theta(u-), \theta, u)\mathcal{P}_\theta(u,t)$ and $R(t,\theta)$ is a mean zero Gaussian process. Its covariance function is given in Section 3.*
(iv) *Let $E_P N(t)$ be continuous, and let $\theta_0$ be an arbitrary point in $\Theta$. If $\hat{\theta}$ is a $\sqrt{n}$-consistent estimate of it, then the process $\hat{W}_0 = \{\hat{W}_0(t) : t \leq \tau\}$, $\hat{W}_0 = \sqrt{n}[\Gamma_{n\hat{\theta}} - \Gamma_{\theta_0} - (\hat{\theta} - \theta_0)\dot{\Gamma}_{\hat{\theta}}]$ converges weakly in $\ell^\infty([0,\tau])$ to $W_0 = W(\cdot, \theta_0)$.*

The proof of this proposition can be found in Section 6.

### 2.3. Some auxiliary notation

From now on we assume that the function $EN(t)$ is continuous. We shall need some auxiliary notation. Define

$$e[f](u,\theta) = \frac{E\{Y(u)[f\alpha](\Gamma_\theta(u),\theta,Z)\}}{E\{Y_i(u)\alpha(\Gamma_\theta(u),\theta,Z)\}},$$

where $f(x,\theta,Z)$, is a function of covariates. Likewise, for any two such functions, $f_1$ and $f_2$, let $\text{cov}[f_1,f_2](u,\theta) = e[f_1 f_2^T](u,\theta) - (e[f_1]e[f_2]^T)(u,\theta)$ and $\text{var}[f](u,\theta) = \text{cov}[f,f](u,\theta)$. We shall write

$$\begin{aligned}
e(u,\theta) &= e[\ell'](u,\theta), & \bar{e}(u,\theta) &= e[\dot{\ell}](u,\theta), \\
v(u,\theta) &= \text{var}[\ell'](u,\theta), & \bar{v}(u,\theta) &= \text{var}[\dot{\ell}](u,\theta), & \rho(u,\theta) &= \text{cov}[\dot{\ell},\ell'](u,\theta),
\end{aligned}$$

for short. Further, let

(2.4)
$$\begin{aligned}
K_\theta(t,t') &= \int_0^{t \wedge t'} C_\theta(du)\mathcal{P}_\theta(u,t)\mathcal{P}_\theta(u,t'), \\
B_\theta(t) &= \int_0^t v(u,\theta)EN(du)
\end{aligned}$$

and define

(2.5)
$$\kappa_\theta(\tau) = \int\int_{0 < u \leq t \leq \tau} C_\theta(du)\mathcal{P}_\theta(u,t)^2 B_\theta(dt).$$

This constant is finite for any point $\tau$ satisfying the condition 2.0 (iii), but is in general infinite, if $\tau_0$ is a continuity point of the survival function $EY(t)$. Finally, we set

$$\begin{aligned}
v_\varphi(t,\theta) &= \bar{v}(t,\theta) + v(t,\theta)\varphi_\theta^{\otimes 2}(t) - \rho(t,\theta)\varphi_\theta^T(t) - \varphi_\theta(t)\rho(t,\theta)^T, \\
\rho_\varphi(t,\theta) &= \rho(t,\theta) - v(t,\theta)\varphi_\theta(t),
\end{aligned}$$

for any function $\varphi_\theta$ square integrable with respect to $B_\theta$. Under the added condition 2.2, we have

$$\begin{aligned}
e(u,\theta_0) &= E[\ell'(\Gamma_{\theta_0}(X),\theta_0,Z)|X=u,\delta=1], \\
v(u,\theta_0) &= \text{var}\,[\ell'(\Gamma_{\theta_0}(X),\theta_0),Z|X=u,\delta=1], \\
\bar{e}(u,\theta_0) &= E[\dot{\ell}(\Gamma_{\theta_0}(X),\theta_0,Z)|X=u,\delta=1], \\
\bar{v}(u,\theta_0) &= \text{var}\,[\dot{\ell}(\Gamma_{\theta_0}(X),\theta_0,Z)|X=u,\delta=1], \\
\rho(u,\theta_0) &= \text{cov}[\dot{\ell}(\Gamma_{\theta_0}(X),\theta_0,Z),\ell'(\Gamma_{\theta_0}(X),\theta_0,Z)|X=u,\delta=1].
\end{aligned}$$

Similarly,

$$\begin{aligned}
v_\varphi(u,\theta_0) &= \text{var}[\dot{\ell}(\Gamma_{\theta_0}(X),\theta_0,Z) - \ell'(\Gamma_{\theta_0}(X),\theta_0,Z)\varphi_{\theta_0}(X)|X=u,\delta=1], \\
\rho_\varphi(u,\theta_0) &= \text{cov}[\dot{\ell}(\Gamma_{\theta_0}(X),\theta_0,Z) \\
&\quad - \ell'(\Gamma_{\theta_0}(X)\theta_0,Z)\varphi_{\theta_0}(X),\ell'(\Gamma_{\theta_0}(X),\theta_0,Z)|X=u,\delta=1].
\end{aligned}$$

However, $e[f], \text{var}[f]$ and $\text{cov}[f,g]$ form conditional expectation and variance–covariance operators even when this assumption fails. This observation, the Cauchy-Schwarz inequality, and the monotone convergence theorem can be used to verify the next lemma.

**Lemma 2.1.** *Suppose that Conditions* 2.0 *and* 2.1 *are satisfied. Let* $EN(t)$ *be continuous, and let* $v(u, \theta) \not\equiv 0$ *a.e.* $EN$.

(i) *If* $\kappa_\theta(\tau_0) < \infty$ *then the kernel* $K_\theta$ *is square integrable with respect to* $B_\theta$. *In addition, if* $e[(\dot{\ell})^{\otimes 2}](u, \theta) \in L_1(EN)$ *then* $\dot{\Gamma}_\theta \in L_2(B_\theta)$.

(ii) *Suppose that the integrability conditions of part (i) are satisfied. For any vector valued function* $\varphi_\theta(t) = \int_0^t g_\theta d\Gamma_\theta \in L_2(B_\theta)$ *the matrices*

$$\Sigma_{0,\varphi}(\theta, \tau) = \int_0^\tau v_\varphi(t, \theta) EN(dt),$$

$$\Sigma_{1,\varphi}(\theta, \tau) = \Sigma_{0,\varphi}(\theta, \tau) + \int_0^\tau \rho_\varphi(t, \theta)[\dot{\Gamma}_\theta(t) + \varphi_\theta(t)]^T EN(dt),$$

$$\Sigma_{2,\varphi}(\theta, \tau) = \Sigma_{0,\varphi}(\theta, \tau)$$
$$+ \int_0^\tau \int_0^\tau K_\theta(t, u) \rho_\varphi(t, \theta) \rho_\varphi(u, \theta)^T EN(du) EN(dt)$$

*have finite components for any point* $\tau \leq \tau_0$.

Here $L_1(EN)$ is the space of functions integrable with respect to $EN$ and $L_2(B_\theta)$ is the space of functions square integrable with respect to $B_\theta$.

**Remark 2.1.** For $\varphi_\theta = -\dot{\Gamma}_\theta$, we have $\Sigma_{1,\varphi}(\theta, \tau) = \Sigma_{2,\varphi}(\theta, \tau)$ if $\rho_{-\dot{\Gamma}}(u, \theta) \equiv 0$ and $v(u, \theta) \not\equiv 0$ a.e. $EN$. If $v(u, \theta) \equiv 0$, then for the sake of completeness, we define $\Sigma_{1,\varphi}(\theta, \tau) = \Sigma_{2,\varphi}(\theta, \tau) = \int_0^\tau \bar{v}(u, \theta) EN(du)$. In this case $\ell_i'(x, \theta, z)$ is a function not depending on covariates. In particular, in the proportional hazards model, we have, $\ell_i'(x, \theta, z) \equiv 0$ for all $\theta$. In scale regression models with hazards $\alpha(x, \theta, z) = \exp(\theta^T z) \alpha_0(x \exp(\theta^T z))$, where $\alpha_0$ is a known function, we have $\ell_i'(x, \theta, z) = \alpha_0'(x)/\alpha_0(x)$ for $\theta = 0$ (independence).

We shall assume now that the point $\tau$ satisfies the condition 2.0 (iii). With this choice, any function $\varphi_\theta(t) = \int_0^t g_\theta d\Gamma_\theta$ of bounded variation on $[0, \tau]$ is square integrable with respect to the measure $B_\theta$, restricted to the interval $[0, \tau]$. However, in Proposition 2.3, we shall allow also for $\tau = \tau_0$ to be a continuity point of the survival function $EY_P(t)$ and assume integrability conditions of Lemma 2.1.

### *2.4. Estimation of the Euclidean component of the model*

To estimate the parameter $\theta$, we use a solution to the score equation $U_n(\theta) = U_{n\varphi_n}(\theta) = 0$, where

$$(2.6) \quad U_{n\varphi_n}(\theta) = \frac{1}{n} \sum_{i=1}^n \int_0^\tau [b_{1i}(\Gamma_{n\theta}(t), t, \theta) - b_{2i}(\Gamma_{n\theta}(t), t, \theta) \varphi_{n\theta}(t)] N_i(dt),$$

$$b_{1i}(x, t, \theta) = \dot{\ell}(x, \theta, Z_i) - [\dot{S}/S](x, \theta, t),$$
$$b_{2i}(x, t, \theta) = \ell'(x, \theta, Z_i) - [S'/S](x, \theta, t)$$

and $\varphi_{n\theta}(t)$ is an estimate of a function $\varphi_\theta(t) = \int_0^t g_\theta d\Gamma_\theta$. We shall make the following regularity assumption.

**Condition 2.3.** Suppose that Conditions 2.0–2.2 hold, and let $\|\cdot\|_v$ be the variation norm on the interval $[0, \tau]$. Let $B(\theta_0, \varepsilon_n) = \{\theta : |\theta - \theta_0| \leq \varepsilon_n\}$ for some sequence $\varepsilon_n \downarrow 0, \sqrt{n} \varepsilon_n \to \infty$. In addition

(i) The matrix $\Sigma_{0,\varphi}(\theta_0, \tau)$ is positive definite.
(ii) The matrix $\Sigma_{1,\varphi}(\theta_0, \tau)$ is non-singular.
(iii) The function $\varphi_{\theta_0}(t) = \int_0^t g_{\theta_0} d\Gamma_{\theta_0}$ satisfies $\|\varphi_{\theta_0}\|_v = O(1)$,
(iv) $\|\varphi_{n\theta_0} - \varphi_{\theta_0}\|_\infty \to_P 0$ and $\limsup_n \|\varphi_{n\theta_0}\|_v = O_P(1)$.
(v) We have either

  (v.1) $\varphi_{n\theta} - \varphi_{n\theta'} = (\theta - \theta')\psi_{n\theta,\theta'}$, where
  $\limsup_n \sup\{\|\psi_{n\theta,\theta'}\|_v : \theta, \theta' \in B(\theta_0, \varepsilon_n)\} = O_P(1)$ or
  
  (v.2) $\limsup_n \sup\{\|\varphi_{n\theta}\|_v : \theta \in B(\theta_0, \varepsilon_n)\} = O_P(1)$ and
  $\sup\{\|\varphi_{n\theta} - \varphi_{\theta_0}\|_\infty : \theta \in B(\theta_0, \varepsilon_n)\} = o_P(1)$.

**Proposition 2.2.** *Suppose that Conditions* 2.3(i)–(iv) *hold.*

(i) *For any $\sqrt{n}$ consistent estimate $\hat{\theta}$ of the parameter $\theta_0$, $\hat{W}_0 = \sqrt{n}[\Gamma_{n\hat{\theta}} - \Gamma_{\theta_0} - (\hat{\theta} - \theta_0)\dot{\Gamma}_{\hat{\theta}}]$ converges weakly in $\ell^\infty([0, \tau])$ to a mean zero Gaussian process $W_0$ with covariance function $\mathrm{cov}(W_0(t), W_0(t')) = K_{\theta_0}(t, t')$.*

(ii) *Suppose that Condition* 2.3(v.1) *is satisfied. Then, with probability tending to 1, the score equation $U_{n\varphi_n}(\theta) = 0$ has a unique solution $\hat{\theta}$ in $B(\theta_0, \varepsilon_n)$. Under Condition* 2.3(v.2), *the score equation $U_{n\varphi_n}(\theta) = o_P(n^{-1/2})$ has a solution, with probability tending to 1.*

(iii) *Define $[\hat{T}, \hat{W}_0], \hat{T} = \sqrt{n}(\hat{\theta} - \theta_0)$, $\hat{W}_0 = \sqrt{n}[\Gamma_{n\hat{\theta}} - \Gamma_{\theta_0} - (\hat{\theta} - \theta_0)\dot{\Gamma}_{\hat{\theta}}]$, where $\hat{\theta}$ are the estimates of part* (ii). *Then $[\hat{T}, \hat{W}_0]$ converges weakly in $R^p \times \ell^\infty([0, \tau])$ to a mean zero Gaussian process $[T, W_0]$ with covariance $\mathrm{cov}\, T = \Sigma_1^{-1}(\theta_0, \tau)\Sigma_2(\theta_0, \tau)[\Sigma_1^{-1}(\theta_0, \tau)]^T$ and*

$$\mathrm{cov}(T, W_0(t)) = -\Sigma_1^{-1}(\theta_0, \tau) \int_0^\tau K_{\theta_0}(t, u)\rho_\varphi(u, \theta_0) EN(du).$$

*Here the matrices $\Sigma_{q,\varphi}, q = 1, 2$ are defined as in Lemma* 2.2.

(iv) *Let $\tilde{\theta}_0$ be any $\sqrt{n}$ estimate, and let $\hat{\varphi}_n = \varphi_{n\tilde{\theta}_0}$ be an estimator of the function $\varphi_{\theta_0}$ such that $\|\hat{\varphi}_n - \varphi_{\theta_0}\|_\infty = o_P(1)$ and $\limsup_n \|\hat{\varphi}_n\|_v = O_P(1)$. Define a one-step M-estimator $\hat{\theta} = \tilde{\theta}_0 + \Sigma_{1\hat{\varphi}_n}(\tilde{\theta}_0, \tau)^{-1}U_{n\hat{\varphi}_n}(\tilde{\theta}_0)$, where $\Sigma_{1,\hat{\varphi}_n}$ is the plug-in analogue of the matrix $\Sigma_{1,\varphi}(\theta_0, \tau)$. Then part* (iii) *holds for the one-step estimator $\hat{\theta}$.*

The proof of this proposition is postponed to Section 7.

**Example 2.1.** A simple choice of the $\varphi_\theta$ function is provided by $\varphi_\theta \equiv 0 = \varphi_{n\theta}$. The resulting score equation, is approximately equal to

$$\hat{U}_n(\theta) = \frac{1}{n}\sum_{i=1}^n \left[ N_i(\tau)\dot{\ell}(\Gamma_{n\theta}(X_i), \theta, Z_i) - \dot{A}(\Gamma_{n\theta}(X_i \wedge \tau), \theta, Z_i) \right],$$

and this score process may be easier to compute in some circumstances. If the transformation $\Gamma$ had been known, the right-hand side would have represented the MLE score function for estimation of the parameter $\theta$. Using results of section 5, we can show that solving equation $\hat{U}_n(\theta) = 0$ or $\hat{U}_n(\theta) = o_P(n^{-1/2})$ for $\theta$ leads to an M estimator asymptotically equivalent to the one in Proposition 2.2. However, this equivalence holds only at rate $\sqrt{n}$. In particular, at the true $\theta_0$, the two score processes satisfy $\sqrt{n}|\hat{U}_n(\theta_0) - U_n(\theta_0)| = o_P(1)$, but they have a different higher order expansions.

**Example 2.2.** The second possible choice corresponds to $\varphi_\theta = -\dot{\Gamma}_\theta$. The score function $U_n(\theta)$ is in this case approximately equal to the derivative of the pseudo-profile likelihood criterion function considered by Bogdanovicius and Nikulin [9] in the case of generalized proportional hazards intensity models. Using results of section 6, we can show that the sample analogue of the function $\dot{\Gamma}_\theta$ satisfies Conditions 2.3(iv) and 2.3(v).

**Example 2.3.** The logarithmic derivatives of $\ell(x, \theta, Z) = \log \alpha(x, \theta, Z)$ may be difficult to compute in some models, so we can try to replace them by different functions. In particular, suppose that $h(x, \theta, Z)$ is a differentiable function with respect to both arguments and the derivatives satisfy a similar Lipschitz continuity assumption as in condition 2.1. Consider the score process (2.6) with function $\varphi_\theta = 0$ and weights $b_{1i}(x, t, \theta) = h(x, \theta, Z_i) - [S_h/S](x, \theta, t)$ where $S_h(x, \theta, t) = \sum_{i=1}^n Y_i(u)[h\alpha](x, \theta, Z_i)$, and $\varphi_{n\theta} \equiv 0$. For $p = 0$ and $p = 2$, define matrices $\Sigma_{p\varphi}^h$ by replacing the functions $v_\varphi$ and $\rho_\varphi$ appearing in matrices $\Sigma_{0\varphi}$ and $\Sigma_{2\varphi}$ with

$$v_\varphi^h(t, \theta_0) = \operatorname{var}[h(\Gamma_{\theta_0}(X), \theta_0, Z)|X = t, \delta = 1],$$
$$\rho_\varphi^h(t, \theta_0) = \operatorname{cov}[h(\Gamma_{\theta_0}(X_i), \theta_0, Z_i), \ell'(\Gamma_{\theta_0}(X_i), \theta_0, Z_i)|X = t, \delta = 1].$$

The matrix $\Sigma_{1\varphi}(\theta_0, \tau)$ is changed to $\Sigma_{1\varphi}^h(\theta_0, \tau) = \int \bar{\rho}_\varphi^h(t, \theta_0) EN(du)$, where the integrand is equal to

$$\operatorname{cov}[h(\Gamma_{\theta_0}(X), \theta_0, Z), \dot{\ell}(\Gamma_{\theta_0}(X), \theta_0, Z) + \ell'(\Gamma_{\theta_0}(X), \theta_0, Z)\dot{\Gamma}_{\theta_0}(X)|X = t, \delta = 1].$$

The statement of Proposition 2.2 remains valid with matrices $\Sigma_{p\varphi}$ replaced by $\Sigma_{p\varphi}^h$, $p = 1, 2$, provided in Condition 2.3 we assume that the matrix $\Sigma_{0\varphi}^h$ is positive definite and the matrix $\Sigma_{1\varphi}^h$ is non-singular. The resulting estimates have a structure analogous to that of the M-estimates considered in the case of uncensored data by Bickel *et al.* [6] and Cuzick [13]. Alternatively, instead of functions $\dot{\ell}_i(x, \theta, z)$ and $\ell'(x, \theta, z)$, the weight functions $b_{1i}$ and $b_{2i}$ can use logarithmic derivatives of a different distribution, with the same parameter $\theta$. The asymptotic variance is of similar form as above. In both cases, the derivations are similar to Section 7, so we do not consider analysis of these score processes in any detail.

**Example 2.4.** Our final example shows that we can choose the $\varphi_\theta$ function so that the asymptotic variance of the estimate $\hat{\theta}$ is equal to the inverse of the asymptotic variance of the normalized score process, $\sqrt{n}U_n(\theta_0)$. Remark 2.1 implies that if $\rho_{-\dot{\Gamma}}(u, \theta_0) \equiv 0$ but $v(u, \theta_0) \not\equiv 0$ a.e. $EN$, then for $\varphi_\theta = -\dot{\Gamma}_\theta$ the matrices $\Sigma_{q,\varphi}$, $q = 1, 2$ are equal. This also holds for $v(u, \theta_0) \equiv 0$. We shall consider now the case $v(u, \theta_0) \not\equiv 0$ and $\rho_{-\dot{\Gamma}}(u, \theta_0) \not\equiv 0$ a.e. $EN$, and without loss of generality, we shall assume that the parameter $\theta$ is one dimensional.

We shall show below that the equation

$$\begin{aligned}(2.7)\quad &\varphi_\theta(t) + \int_0^\tau K_\theta(t, u)v(u, \theta)\varphi_\theta(u)EN(du) \\ &= -\dot{\Gamma}_\theta(t) + \int_0^\tau K_\theta(t, u)\rho(u, \theta)EN(du)\end{aligned}$$

has a unique solution $\varphi_\theta$ square integrable with respect to the measure (2.4). For $\theta = \theta_0$, the corresponding matrices $\Sigma_{1,\varphi}(\theta_0, \tau)$ and $\Sigma_{2,\varphi}(\theta_0, \tau)$ are finite. Substitution of the conditional correlation function $\rho_\varphi(t, \theta_0) = \rho(t, \theta_0) - \varphi_{\theta_0}(t)v(t, \theta_0)$ into the matrix $\Sigma_{2,\varphi}(\theta_0, \tau)$ shows that they are also equal. (In the multiparameter case, the equation (2.7) is solved for each component of the $\theta$).

Equation (2.7) simplifies if we replace the function $\varphi_\theta$ by $\psi_\theta = \varphi_\theta + \dot{\Gamma}_\theta$. We get

$$(2.8)\quad \psi_\theta(t) - \lambda \int_0^\tau K_\theta(t, u)\psi_\theta(u)B_\theta(du) = \eta_\theta(t),$$

where $\lambda = -1$,
$$\eta_\theta(t) = \int_0^\tau K_\theta(t,u)\rho_{-\dot{\Gamma}}(u,\theta)EN(du),$$

$\rho_{-\dot{\Gamma}}(u,\theta) = v(u,\theta)\dot{\Gamma}_\theta(u)+\rho(u,\theta)$ and $B_\theta$ is given by (2.4). For fixed $\theta$, the kernel $K_\theta$ is symmetric, positive definite and square integrable with respect to $B_\theta$. Therefore it can have only positive eigenvalues. For $\lambda = -1$, the equation has a unique solution given by

$$(2.9) \qquad \psi_\theta(t) = \eta_\theta(u) - \int_0^\tau \Delta_\theta(t,u,-1)\eta_\theta(u)B_\theta(du),$$

where $\Delta_\theta(t,u,\lambda)$ is the resolvent corresponding to the kernel $K_\theta$. By definition, the resolvent satisfies a pair of integral equations

$$\begin{aligned} K_\theta(t,u) &= \Delta_\theta(t,u,\lambda) - \lambda\int_0^\tau \Delta_\theta(t,w,\lambda)B_\theta(dw)K_\theta(w,u) \\ &= \Delta_\theta(t,u,\lambda) - \lambda\int_0^\tau K_\theta(t,w)B_\theta(dw)\Delta_\theta(w,u,\lambda), \end{aligned}$$

where integration is with respect to different variables in the two equations. For $\lambda = -1$ the solution to the equation is given by

$$\begin{aligned} \psi_\theta(t) &= \int_0^\tau K_\theta(t,u)\rho_{-\dot{\Gamma}}(u,\theta)EN(du) \\ &\quad - \int_0^\tau \Delta_\theta(t,w,-1)B_\theta(dw)\int_0^\tau K_\theta(w,u)\rho_{-\dot{\Gamma}}(u,\theta)EN(du) \end{aligned}$$

and the resolvent equations imply that the right-hand side is equal to

$$(2.10) \qquad \psi_\theta(t) = \int_0^\tau \Delta_\theta(t,u,-1)\rho_{-\dot{\Gamma}}(u,\theta)EN(du).$$

For $\theta = \theta_0$, substitution of this expression into the formula for the matrices $\Sigma_{1,\varphi}(\theta_0,\tau)$ and $\Sigma_{2,\varphi}(\theta_0,\tau)$ and application of the resolvent equations yields also

$$\begin{aligned} \Sigma_{1,\varphi}(\theta_0,\tau) &= \Sigma_{2,\varphi}(\theta_0,\tau) \\ &= \int_0^\tau v_{-\dot{\Gamma}}(u,\theta_0)EN(du) \\ &\quad - \int_0^\tau\int_0^\tau \Delta_{\theta_0}(t,u,-1)\rho_{-\dot{\Gamma}}(u,\theta_0)\rho_{-\dot{\Gamma}}(t,\theta_0)^T EN(du)EN(dt). \end{aligned}$$

It remains to find the resolvent $\Delta_\theta$. We shall consider first the case of $\theta = \theta_0$. To simplify algebra, we multiply both sides of the equation (2.8) by $\mathcal{P}_{\theta_0}(0,t)^{-1} = \exp\int_0^t s'(\theta_0,\Gamma_{\theta_0}(u),u)C_{\theta_0}(du)$. For this purpose set

$$\begin{aligned} \tilde{\psi}(t) &= \mathcal{P}_{\theta_0}(0,t)^{-1}\psi(t), & \dot{G}(t) &= \mathcal{P}_{\theta_0}(0,t)^{-1}\dot{\Gamma}_{\theta_0}(t), \\ \tilde{v}(t,\theta_0) &= v(t,\theta_0)\mathcal{P}_{\theta_0}(0,t)^2, & \tilde{\rho}_{-\dot{G}}(t,\theta_0) &= \mathcal{P}_{\theta_0}(0,t)\rho_{-\dot{\Gamma}}(t,\theta_0), \\ b(t) &= \int_0^t \tilde{v}(u,\theta_0)dEN(u), & c(t) &= \int_0^t \mathcal{P}_{\theta_0}(0,u)^{-2}dC_{\theta_0}(u). \end{aligned}$$

Multiplication of (2.8) by $\mathcal{P}_{\theta_0}(0,t)^{-1}$ yields

$$(2.11) \qquad \tilde{\psi}(t) + \int_0^\tau k(t,u)\tilde{\psi}(u)b(du) = \int_0^\tau k(t,u)\tilde{\rho}_{-\dot{G}}(u,\theta_0)EN(du),$$

where the kernel $k$ is given by $k(t, u) = c(t \wedge u)$. Since this is the covariance function of a time transformed Brownian motion, we obtain a simpler equation. The solution to this Fredholm equation is

$$(2.12) \qquad \tilde{\psi}(t) = \int_0^\tau \tilde{\Delta}(t, u) \tilde{\rho}_{-\dot{G}}(u, \theta_0) EN(du),$$

where $\tilde{\Delta}(t, u) = \tilde{\Delta}(t, u, -1)$, and $\tilde{\Delta}(t, u, \lambda)$ is the resolvent corresponding to the kernel $k$. More generally, we consider the equation

$$(2.13) \qquad \tilde{\psi}(t) + \int_0^\tau k(t, u) \tilde{\psi}(u) b(du) = \tilde{\eta}(t).$$

Its solution is of the form

$$\tilde{\psi}(t) = \tilde{\eta}(t) - \int_0^\tau \tilde{\Delta}(t, u) b(du) \tilde{\eta}(u).$$

To give the form of the $\tilde{\Delta}$ function, note that the constant $\kappa_{\theta_0}(\tau)$ defined in (2.5) satisfies

$$\kappa(\tau) = \kappa_{\theta_0}(\tau) = \int_0^\tau c(u) b(du).$$

**Proposition 2.3.** *Suppose that Assumptions 2.0(i) and (ii) are satisfied and $v(u, \theta_0) \not\equiv 0$, For $j = 0, 1, 2, 3$, $n \geq 1$ and $s < t$ define interval functions $\Psi_j(s, t) = \sum_{m=0}^\infty \Psi_{jm}(s, t)$ as follows:*

$$\Psi_{00}(s, t) = 1, \quad \Psi_{20}(s, t) = 1,$$

$$\Psi_{0n}(s, t) = \int\!\!\int_{s < u_1 \leq u_2 \leq t} \Psi_{0,n-1}(s, u_1-) c(du_1) b(du_2) \quad n \geq 1,$$

$$\Psi_{1n}(s, t) = \int_{(s,t]} \Psi_{0n}(s, u-) c(du) \quad n \geq 0,$$

$$\Psi_{2n}(s, t) = \int\!\!\int_{s \leq u_1 < u_2 < t} b(du_1) c(du_2) \Psi_{2,n-1}(u, t)$$

$$= \int_{[s,t)} b(du_1) \Psi_{1,n-1}(u, t-) \quad n \geq 1,$$

$$\Psi_{3n}(s, t) = \int_{[s,t)} \Psi_{2n}(s, u) b(dw) \quad n \geq 0.$$

*For $j = 2, 3$, define $\Psi_{jn}(s, t+)$ by replacing the intervals $[s, t)$ with $[s, t]$ in the last two lines, and similarly, define $\Psi_{jn}(s, t-)$ by replacing intervals $(s, t]$ with $(s, t)$ in the first two definitions. For $s > t$, set $\Psi_j(s, t) = 0, j = 0, 1, 2, 3$ and let $\Psi_{j0}(t, t) = 1$ for $j = 0, 2$, $\Psi_{10}(t, t) = c(\Delta t)$, $\Psi_{30}(t, t) = b(\Delta t)$, and $\Psi_{jn}(t, t) = 0$ for $n \geq 1$, $j = 0, 1, 2, 3$.*

(i) *We have*

$$\Psi_0(s, t) = 1 + \int_{(s,t]} \Psi_1(s, u) b(du) = 1 + \int_{(s,t]} c(du) \Psi_3(u, t+),$$

$$\Psi_1(s, t) = \int_{(s,t]} \Psi_0(s, u-) c(du) = \int_{(s,t]} c(du) \Psi_2(u, t+),$$

$$\Psi_2(s, t) = 1 + \int_{[s,t)} b(du) \Psi_1(u, t-) = 1 + \int_{[s,t)} \Psi_3(s, u) c(du),$$

$$\Psi_3(s, t) = \int_{[s,t)} \Psi_2(s, u) b(du) = \int_{[s,t)} b(du) \Psi_0(u, t-).$$

*For any point $\tau$ satisfying Condition 2.0(iii), $\Psi_j, j = 0, 1, 2, 3$ form bounded monotone increasing interval functions. In particular, $\Psi_0(s,t) \leq \exp\kappa(\tau)$ and $\Psi_1(s,t) \leq \Psi_0(s,t)[c(t) - c(s)]$. In addition if $\tau_0$ is a continuity point of the survival function $E_P Y(t)$ and $\kappa(\tau_0) < \infty$, then $\Psi_0(s,t) \leq \exp\kappa(\tau_0)$ for any $0 < s < t \leq \tau_0$, while the remaining functions are locally bounded.*

(ii) *Suppose that $\tau$ satisfies Condition 2.0(iii), or else $\tau = \tau_0$, $\tau_0$ is a continuity point of $E_P Y(t)$ and $\kappa(\tau_0) < \infty$. The resolvent of the kernel $k$ is given by $\tilde{\Delta}(s,t,-1) = \tilde{\Delta}(s,t) = \Psi_0(0,\tau)^{-1}\Psi_1(0, s \wedge t)\Psi_0(s \vee t, \tau)$.*

(iii) *Under assumptions in (ii), for any $\tilde{\eta} \in L_2(b)$, the solution to equation (2.13) satisfies $\tilde{\psi} \in L_2(b)$, and and $\|\tilde{\psi}\|_2 \leq \|\tilde{\eta}\|_2[1 + \Psi_0(0,\tau_0)\kappa(\tau_0)]$, where $\|\cdot\|_2$ is the $L_2$ norm with respect to the measure $b$.*

(iv) *Suppose that $\tau$ satisfies Condition 2.0(iii). If $\tilde{\eta}$ is a bounded function or a function of bounded variation, then the solution $\tilde{\psi}$ has the same properties and the bounds of part (ii) hold in supremum and variation norm, respectively.*

(v) *The solution to equation (2.7) is given by*

$$\varphi_{\theta_0}(t) = -\dot{\Gamma}_{\theta_0}(t) + \int_0^\tau \tilde{\Delta}(t,u)\rho_{-\dot{\Gamma}}(u,\theta_0)EN(du)\mathcal{P}_{\theta_0}(0,u)\mathcal{P}_{\theta_0}(0,t).$$

*Under assumptions of part (ii) and integrability conditions of Lemma 2.2, we have $\varphi \in L_2(B_{\theta_0})$. We also have, $\varphi_{\theta_0}(t) = \int_0^t g_{\theta_0} d\Gamma_{\theta_0}$, where*

$$g_{\theta_0}(u) = [\dot{s}/s](\Gamma_{\theta_0}(u), \theta_0, u) - [s'/s](\Gamma_{\theta_0}(u), \theta_0, u)\varphi_{\theta_0}(u) - s(\Gamma_{\theta_0}(u), \theta_0, u)^{-1}\int_u^\tau \mathcal{P}_{\theta_0}(u,t)\rho_\varphi(u,\theta_0)EN(du),$$

(vi) *If $\tau$ satisfies Condition 2.0(iii), then the solution $\varphi$ is a function of bounded variation. Moreover, the constant $W = \Psi(0,\tau)$ satisfies $W = \Psi_1(0,t) \times \Psi_3(t,\tau) + \Psi_0(0,t)\Psi_0(t,\tau)$ for any $0 < t \leq \tau$, and*

$$\varphi_{\theta_0}(t) = \int_0^\tau \tilde{\Delta}(t,u)\rho(u,\theta_0)EN(du)\mathcal{P}_{\theta_0}(0,u)\mathcal{P}_{\theta_0}(0,t) + \int_0^\tau \bar{\Delta}(t,u)\dot{s}(\Gamma_0(u),\theta_0,u)c(du)\mathcal{P}_{\theta_0}(0,u)\mathcal{P}_{\theta_0}(0,t),$$

*where $\bar{\Delta}(t,u) = W^{-1}[\Psi_0(0, u \wedge t)\Psi_0(u \vee t, \tau) - \Psi_1(0, u \wedge t)\Psi_3(u \vee t, \tau)]$.*

The proof of this proposition is given in Section 8.

We have chosen to transform equation (2.8) in order to simplify calculations. The resolvent of the kernel $K$ corresponding to equation (2.8) can be obtained based on recurrent Fredholm determinant formulas [25] applied to the kernel $K$. The same arguments can be applied to find the solution to equation (2.8) for $\theta \neq \theta_0$. The only difference is that the kernel function $K_\theta(t,u)$ does not represent the asymptotic covariance function of the process $\sqrt{n}[\Gamma_{n\theta} - \Gamma_\theta]$ for such $\theta$ points.

The sample analogue of the function $\psi_\theta$ can be obtained in several different manners. Firstly, equations (2.7)–(2.8) can be solved directly by plugging in sample analogues of the functions $K, v, \rho$ etc. If these sample analogues are functions placing mass at each uncensored observation, then this choice is not convenient, because to solve the equation one must eventually invert an $m \times m$ dimensional matrix (here $m$ is the number of uncensored observations in the sample). Proposition 2.3 provides a simpler form of this equation. Define estimates

$$c_{n\theta}(t) = \int_0^t \tilde{\mathcal{P}}_{n\theta}(0,u)^{-2} C_{n\theta}(du),$$

$$b_{n\theta}(t) = \int_0^t \tilde{\mathcal{P}}_{n\theta}(0,u)^2 B_{n\theta}(du),$$

$$C_{n\theta}(du) = \int_0^t S(\Gamma_{n\theta}(u-),\theta,u)^{-2} N_{\cdot}(du),$$

$$\tilde{\mathcal{P}}_{n\theta}(u,t) = \exp[-\int_u^t S'(\Gamma_\theta(u-),\theta,u) C_{n\theta}(du)],$$

and let $B_{n\theta}$ be the plug-in analogue of the formula (2.4). Let $X^*_{(1)} < \cdots < X^*_{(m)}$ be the distinct ordered uncensored observations in the sample. Then the discrete version of equations (2.11)–(2.13) is given by

$$\tilde{\psi}_{n\theta}(X^*_{(j)}) + \sum_{i=1}^m c_{n\theta}(X^*_{(i)} \wedge X^*_{(j)}) b_{n\theta}(\Delta X^*_{(i)}) \tilde{\psi}_{n\theta}(X^*_{(i)}) = \tilde{\eta}_{n\theta}(X^*_{(j)}).$$

Using Proposition 2.3, we have shown in an earlier version of this text that finding solution to this equation amounts to inversion of a bandsymmetric tridiagonal matrix which can be easily implemented in practice. A numerical example is given in [16]. The formula (2.14) in part (v) gives in this case an estimate of the function $\varphi_\theta$ corresponding to the equation (2.7). We show in [17] that it satisfies Conditions 2.3.

Finally, we show that in the continuous case, equation (2.11) corresponds to a Sturm–Liouville equation. Suppose that the point $\tau$ satisfies Condition 2.0(iii). By twice "differentiating" (2.11) with respect to $dc$, we obtain a Sturm–Liouville equation

$$\frac{d}{dc}[\frac{d}{dc}\tilde{\psi}](t) = \tilde{\psi}(t)\frac{db}{dc}(t) - \tilde{\rho}_{-\dot{G}}(t,\theta_0)\frac{EN}{dc}(t)$$

with boundary conditions $\tilde{\psi}(0) = 0, \frac{d}{dc}\tilde{\psi}(t)_{|t=\tau} = 0$. Its solution is of the form (2.12) with $\tilde{\Delta}$ representing the Green's function associated with the homogeneous equation

$$(2.14) \qquad \frac{d}{dc}[\frac{d}{dc}\tilde{\psi}](t) = \tilde{\psi}(t)\frac{db}{dc}(t)$$

and boundary conditions $\tilde{\psi}(0) = 0, \frac{d}{dc}\tilde{\psi}(\tau) = 0$. The Green's function is given by $\tilde{\Delta}(t,u) = \tilde{W}^{-1}[\psi_1(t \wedge u)\psi_0(t \vee u)]$, where $\psi_1$ and $\psi_0$ is a pair of fundamental solutions, $\psi_1$ corresponding to the left boundary ($\psi_1(0) = 0$) and $\psi_0$ corresponding to the right-boundary ($\frac{d}{dc}\psi_0(\tau) = 0$). Moreover,

$$\tilde{W} = -[\psi_1(t)\frac{d}{dc}\psi_0(t) - \psi_0(t)\frac{d}{dc}\psi_1(t)]$$

is the negative Wronskian (the right-hand side is a constant, not depending on t). By twice integrating the homogeneous equation subject to the boundary conditions, we obtain a pair of Volterra equations whose solutions are $a_1\Psi_1(0,t)$ and $a_0\Psi_0(t,\tau)$, where $a_p \neq 0$ are arbitrary constants. The choice of $a_p = 1$, corresponds to the Volterra equations for $\Psi_0(t,\tau)$ and $\Psi_1(0,t)$ discussed in part (i). We also have $\tilde{W} = a_0a_1\Phi_0(0,\tau) = a_0a_1W$. Thus the $\tilde{\Delta}$ function of Proposition 2.3 is the Green's function of this Sturm–Liouville equation. Note that this is a different equation than in Bickel [4] and Bickel *et al.* [6]. In particular, it derives its form from the covariance function of a time transformed Brownian motion, rather than Brownian Bridge.

## 3. Examples

In this section we assume the conditional independence Assumption 2.2 and discuss Condition 2.1(ii) in more detail. It assumes that the hazard rate satisfies $m_1 \leq \alpha(x,\theta,z) \leq m_2$ $\mu$ a.e. $z$. This holds for example in the proportional hazards model, if the covariates are bounded and the regression coefficients vary over a bounded neighbourhood of the true parameter. Recalling that for any $P \in \mathcal{P}$, $\mathcal{X}(\mathcal{P})$ is the set of (sub)-distribution functions whose cumulative hazards satisfy $m_2^{-1}A \leq g \leq m_1^{-1}A$ and $A$ is the cumulative hazard function (2.1), this uniform boundedness is used in Section 6 to verify that equation (2.2) has a unique solution which is defined on the entire support of the withdrawal time distribution. This need not be the case in general, as the equation may have an explosive solution on an interval strictly contained in the support of this distribution ([20]).

We shall consider now the case of hazards $\alpha(x,\theta,z)$ which for $\mu$ almost all $z$ are locally bounded and locally bounded away from 0. A continuous nonnegative function $f$ on the positive half-line is referred to here as locally bounded and locally bounded away from 0, if $f(0) > 0$, $\lim_{x \uparrow \infty} f(x)$ exists, and for any $d > 0$ there exists a finite positive constant $k = k(d)$ such that $k^{-1} \leq f(x) \leq k$ for $x \in [0,d]$. In particular, hazards of this form may form unbounded functions growing to infinity or functions decaying to 0 as $x \uparrow \infty$.

To allow for this type of hazards, we note that the transformation model assumes only that the conditional cumulative hazard function of the failure time $T$ is of the form $H(t|z) = A(\tilde{\Gamma}(t),\theta|z)$ for some *unspecified* increasing function $\tilde{\Gamma}$. We can choose it as $\tilde{\Gamma} = \Phi(\Gamma)$, where $\Phi$ is a known increasing differentiable function mapping positive half-line onto itself, $\Phi(0) = 0$. This is equivalent to selection of the reparametrized core model with cumulative hazard function $A(\Phi(x),\theta|z)$ and hazard rate $\alpha(\Phi(x),\theta|z)\varphi(\Phi(x))$, $\varphi = \Phi'$. If in the original model the hazard rate decays to 0 or increases to infinity at its tails, then in the reparametrized model the hazard rate may form a bounded function. Our results imply in this case that we can define a family of transformations $\tilde{\Gamma}_\theta$ bounded between $m_2^{-1}A(t)$ and $m_1^{-1}A(t)$, This in turn defines a family of transformations $\Gamma_\theta$ bounded between $\Phi^{-1}(m_2^{-1}A(t))$ and $\Phi^{-1}(m_1^{-1}A(t))$. More generally, the function $\Phi$ may depend on the unknown parameter $\theta$ and covariates. Of course selection of this reparametrization is not unique, but this merely means that different core models may generate the same semiparametric transformation model.

**Example 3.1.** Half-logistic and half-normal scale regression model. The assumption that the conditional distribution of a failure time $T$ given a covariate $Z$ has cumulative hazard function $H(t|z) = A_0(\tilde{\Gamma}(t)\exp[\theta^T z])$, for some unknown increasing function $\tilde{\Gamma}$ (model I), is clearly equivalent to the assumption that this cumulative hazard function is of the form $H(t|z) = A_0(A_0^{-1}(\Gamma(t))\exp[\theta^T z])$, for some unknown increasing function $\Gamma$ (model II). The corresponding core models have hazard rates

$$(3.1) \qquad \text{model I:} \qquad \alpha(x,\theta,z) = e^{\theta^T z}\alpha_0(xe^{\theta^T z})$$

and

$$(3.2) \qquad \text{model II:} \quad \alpha(x,\theta,z) = e^{\theta^T z}\frac{\alpha_0(A_0^{-1}(x)e^{\theta^T z})}{\alpha_0(A_0^{-1}(x))},$$

respectively. In the case of the core model I, Condition 2.1(ii) is satisfied if the covariates are bounded, $\theta$ varies over a bounded neighbourhood of the true parameter and $\alpha_0$ is a hazard rate that is bounded and bounded away from 0. An example is provided by the half-logistic transformation model with $\alpha_0(x) = 1/2 + \tanh(x/2)$. This is a bounded increasing function from $1/2$ to $1$.

Next let us consider the choice of the half-normal transformation model. The half-normal distribution has survival function $F_0(x) = 2(1 - \Phi(x))$, where $\Phi$ is the standard normal distribution function. The hazard rate is given by

$$\alpha_0(x) = x + \frac{\int_x^\infty F_0(u)du}{F_0(x)}.$$

The second term represents the residual mean of the half normal distribution, and we have $\alpha_0(x) = x + \ell'_0(x)$. The function $\alpha_0$ is increasing and unbounded so that the Condition 2.1(ii) fails to be satisfied by hazard rates (3.1). On the other hand the reparameterized transformation model II has hazard rates

$$\alpha(x, \theta, z) = e^{\theta^T z} \frac{A_0^{-1}(x)e^{\theta^T z} + \ell'_0(e^{\theta^T z} A_0^{-1}(x))}{A_0^{-1}(x) + \ell'_0(A_0^{-1}(x))}.$$

It can be shown that the right side satisfies $\exp(\theta^T z) \le \alpha(x, \theta, z) \le \exp(2\theta^T z) + \exp(\theta^T z)$ for $\exp(\theta^T z) > 1$, and $\exp(2\theta^T z)(1+\exp(\theta^T z))^{-1} \le \alpha(x, \theta, z) \le \exp(\theta^T z)$ for $\exp(\theta^T z) \le 1$. These inequalities are used to verify that the hazard rates of the core model II satisfy the remaining conditions 2.1 (ii).

Condition 2.1 assumes that the support of the distribution of the core model corresponds to the whole positive half-line and thus it has a support independent of the unknown parameter. The next example deals with the situation in which this support may depend on the unknown parameter.

**Example 3.2.** The gamma frailty model [14, 28] has cumulative hazard function

$$\begin{aligned}
G(x, \theta|z) &= \frac{1}{\eta} \log[1 + \eta x e^{\beta^T z}], \quad \theta = (\eta, \beta), \eta > 0, \\
&= x e^{\beta^T z}, \quad \eta = 0, \\
&= \frac{1}{\eta} \log[1 + \eta x e^{\beta^T z}], \quad \text{for} \quad \eta < 0 \quad \text{and} \quad -1 < \eta e^{\beta^T z} x \le 0.
\end{aligned}$$

The right-hand side can be recognized as inverse cumulative hazard rate of Gompertz distribution.

For $\eta < 0$ the model is not invariant with respect to the group of strictly increasing transformations of $R_+$ onto itself. The unknown transformation $\Gamma$ must satisfy the constraint $-1 < \eta \exp(\beta^T z) \Gamma(t) \le 0$ for $\mu$ a.e. $z$. Thus its range is bounded and depends on $(\eta, \beta)$ and the covariates. Clearly, in this case the transformation model, assuming that the function $\Gamma$ does not depend on covariates and parameters does not make any sense. When specialized to the transformation $\Gamma(t) = t$, the model is also not regular. For example, for $\eta = -1$ the cumulative hazard function is the same as that of the uniform distribution on the interval $[0, \exp(-\beta^T Z)]$. Similarly to the uniform distribution without covariates, the rate of convergence of the estimates of the regression coefficient is $n$ rather than $\sqrt{n}$. For other choices of the $\tilde{\eta} = -\eta$ parameter, the Hellinger distance between densities corresponding to parameters $\beta_1$ and $\beta_2$ is determined by the magnitude of

$$E_Z 1(h^T Z > 0)[1 - \tilde{\eta} \exp(-h^T Z)]^{1/\tilde{\eta}} + E_Z 1(h^T Z < 0)[1 - \tilde{\eta} \exp(h^T Z)]^{1/\tilde{\eta}},$$

where $h = \beta_2 - \beta_1$. After expanding the exponents, this difference is of order $O(E_Z |hZ|^{1/\tilde{\eta}})$ so that for $\tilde{\eta} \le 1/2$ the model is regular, and irregular for $\tilde{\eta} > 1/2$.

For $\eta \geq 0$, the model is Hellinger differentiable both in the presence of covariates and in the absence of them ($\beta = 0$). The densities are supported on the whole positive half-line. The hazard rates are given by $g(x, \theta|z) = \exp(\beta^T z)\,[1 + \eta \exp(\beta^T z)x]^{-1}$. These are decreasing functions decaying to zero as $x \uparrow \infty$. Using Gompertz cumulative hazard function $G_\eta^{-1}(x) = \eta^{-1}[e^{\eta x} - 1]$ to reparametrize the model, we get $A(x, \theta|z) = G(G_\eta^{-1}(x), \theta|z) = \eta^{-1}\log[1 + (e^{\eta x} - 1)\exp(\beta^T z)]$. The hazard rate of this model is given by $\alpha(x, \theta|z) = \exp(\beta^T z + \eta x)[1 + (e^{\eta x} - 1)\exp(\beta^T z)]^{-1}$. Pointwise in $\beta$, this function is bounded between $\max\{\exp(e\beta^T Z), 1\}$ and from below by $\min\{\exp(\beta^T Z), 1\}$. The bounds are uniform for all $\eta \in [0, \infty)$ and the reparametrization preserves regularity of the model.

Note that the original core model has the property that for each parameter $\eta, \eta \geq 0$ it describes a distribution with different shape and upper tail behaviour. As a result of this, in the case of transformation model, the unknown function $\Gamma$ is confounded by the parameter $\eta$. For example, at $\eta = 0$, the unknown transformation $\Gamma$ represents a cumulative hazard function whereas at $\eta = 1$, it represents an odds ratio function. For any continuous variable $X$ having a nondefective distribution, we have $E\Gamma(X) = 1$, if $\Gamma$ is a cumulative hazard function, and $E\Gamma(X) = \infty$, if $\Gamma$ is an odds ratio function. Since an odds ratio function diverges to infinity at a much faster rate than a cumulative hazard function, these are clearly very different parameters.

The preceding entails that when $\eta, \eta \geq 0$, is unknown we are led to a constrained optimization problem and our results fail to apply. Since the parameter $\eta$ controls the shape and growth-rate of the transformation, it is not clear why this parameter could be identifiable based on rank statistics instead of order statistics. But if omission of constraints is permissible, then results of the previous section apply so long as the true regression coefficient satisfies $\beta_0 \neq 0$ and there exists a preliminary $\sqrt{n}$-consistent estimator of $\theta$. At $\beta_0 = 0$, the parameter $\eta$ is not identifiable based on ranks, if the unknown transformation is only assumed to be continuous and completely specified. We do not know if such initial estimators exist, and rank invariance arguments used in [14] suggest that the parameter $\eta$ is not identifiable based on rank statistics because the models assuming that the cumulative hazard function is of the form $\eta^{-1}\log[1 + c\eta \exp(\beta^T z)\Gamma(t)]$ and $\eta^{-1}\log[1 + \exp(\beta^T z)\Gamma(t)], c > 0, \eta > 0$ all represent the same transformation model corresponding to log-Burr core model with different scale parameter c. Because this scale parameter is not identifiable based on ranks, the restriction $c = 1$ does not imply, that $\eta$ may be identifiable based on rank statistics.

The difficulties arising in analysis of the gamma frailty with fixed frailty parameter disappear if we assume that the frailty parameter $\eta$ depends on covariates. One possible choice corresponds to the assumption that the frailty parameter is of the form $\eta(z) = \exp \xi^T z$. The corresponding cumulative hazard function is given by $\exp[-\xi^T z]\log[1 + \exp(\xi^T z + \beta^T z)\Gamma(t)]$. This is a frailty model assuming that conditionally on $Z$ and an unobserved frailty variable $U$, the failure time $T$ follows a proportional hazards model with cumulative hazard function $U\Gamma(t)\exp(\beta^T Z)$, and conditionally on $Z$, the frailty variable $U$ has gamma distribution with shape and scale parameter equal to $\exp(\xi^T z)$.

**Example 3.3.** Linear hazard model. The core model has hazard rate $h(x, \theta|z) = a_\theta(z) + xb_\theta(z)$ where $a_\theta(z), b_\theta(z)$ are nonnegative functions of the covariates dependent on a Euclidean parameter $\theta$. The cumulative hazard function is equal to $H(t|z) = a_\theta(z)t + b_\theta(z)t^2/2$. Note that the shape of the density depends on the parameters $a$ and $b$: it may correspond to both a decreasing and a non-monotone function.

Suppose that $b_\theta(z) > 0, a_\theta(z) > 0$. To reparametrize the model we use $G^{-1}(x) = [(1+2x)^{1/2}-1]$. The reparametrized model has cumulative hazard function $A(x,\theta|z) = H(G^{-1}(x), \theta|z)$ with hazard rate $\alpha(x,\theta,z) = a_\theta(z)(1 + 2x)^{-1/2} + b_\theta(z)[1 - (1 + 2x)^{-1/2}]$. The hazard rates are decreasing in $x$ if $a_\theta(z) > b_\theta(z)$, constant in $x$ if $a_\theta(z) = b_\theta(z)$ and bounded increasing if $a_\theta(z) < b_\theta(z)$. Pointwise in $z$ the hazard rates are bounded from above by $\max\{a_\theta(z), b_\theta(z)\}$ and from below by $\min\{a_\theta(z), b_\theta(z)\}$. Thus our regularity conditions are satisfied, so long as in some neighbourhood of the true parameter $\theta_0$ these maxima and minima stay bounded and bounded away from 0 and the functions $a_\theta, b_\theta$ satisfy appropriate differentiability conditions. Finally, a sufficient condition for identifiability of parameters is that at a known reference point $z_0$ in the support of covariates, we have $a_\theta(z_0) = 1 = b_\theta(z_0)$, $\theta \in \Theta$ and

$$[a_\theta(z) = a_{\theta'}(z) \quad \text{and} \quad b_\theta(z) = b_{\theta'}(z) \quad \mu \text{ a.e. } z] \Rightarrow \theta = \theta'.$$

Returning to the original linear hazard model, we have excluded the boundary region $a_\theta(z) = 0$ or $b_\theta(z) = 0$. These boundary regions lead to lack of identifiability. For example,

$$\begin{aligned}
\text{model 1:} &\quad a_\theta(z) = 0 \quad \mu \text{ a.e. } z, \\
\text{model 2:} &\quad b_\theta(z) = 0 \quad \mu \text{ a.e. } z, \\
\text{model 3:} &\quad a_\theta(z) = c b_\theta(z) \quad \mu \text{ a.e. } z,
\end{aligned}$$

where $c > 0$ is an arbitrary constant, represent the same proportional hazards model. The reparametrized model does not include the first two models, but, depending on the choice of the parameter $\theta$, it may include the third model (with $c = 1$).

**Example 3.4.** Half-t and polynomial scale regression models. In this example we assume that the core model has cumulative hazard $A_0(x \exp[\theta^T z])$ for some known function $A_0$ with hazard rate $\alpha_0$. Suppose that $c_1 \le \exp(\theta^T z) \le c_2$ for $\mu$ a.e. $z$.

For fixed $\xi \ge -1$ and $\eta \ge 0$, let $G^{-1}$ be the inverse cumulative hazard function corresponding to the hazard rate $g(x) = [1 + \eta x]^\xi$. If $\alpha_0/g$ is a function locally bounded and locally bounded away from zero such that $\lim_{x \uparrow \infty} \alpha_0(x)/g(x) = c$ for a finite positive constant $c$, then for any $\varepsilon \in (0,c)$ there exist constants $0 < m_1(\varepsilon) < m_2(\varepsilon) < \infty$, such that the hazard rate of $A_0(G^{-1}(x) \exp[\theta^T z])$ is bounded between $m_1(\varepsilon)$ and $m_2(\varepsilon)$. Indeed, using $c_1 \le \exp(\theta^T z) \le c_2$ and monotonicity properties of the function $g(x)$, we can find finite positive constants $b_1, b_2$ such that $b_1 \le e^{\theta^T z} g(x \exp[\theta^T z])/g(x) \le b_2$ for $\mu$ a.e. $z$ and $x \ge 0$. The claim follows by setting $m_1(\varepsilon) = b_1 \max(c-\varepsilon, k^{-1})$ and $m_2(\varepsilon) = b_2 \min(c+\varepsilon, k)$, where $k = k(d), k > 0$ and $d > 0$ are such that $c - \varepsilon \le \alpha_0(x)/g(x) \le c + \varepsilon$ for $x > d$, and $k^{-1} \le \alpha_0(x)/g(x) \le k$, for $x \le d$.

In the case of half-logistic distribution, we choose $g(x) \equiv 1$. The function $g(x) = 1 + x$ applies to the half-normal scale regression, while the choice $g(x) = (1 + n^{-1} x)^{-1}$ applies to the half-$t_n$ scale regression model. Of course in the case of gamma, inverse Gaussian frailty models (with fixed frailty parameters) and linear hazard model the choice of the $g(x)$ function is obvious.

In the case of polynomial hazards $\alpha_0(x) = 1 + \sum_{p=1}^m a_p x^p, m > 1$, where $a_p$ are fixed nonnegative coefficients and $a_m > 0$, we choose $g(x) = [1 + a_m x]^m$. Note however, that polynomial hazards may be also well defined when some of the coefficients $a_p$ are negative. We do not know under what conditions polynomial hazards define regular parametric models, but we expect that in such models parameters are estimated subject to added constraints in both parametric and semiparametric setting. Evidently, our results do not apply to such complicated problems.

The choice of $g(x) = [1 + \eta x]^\xi, \xi < -1$ was excluded in this example because it forms a defective hazard rate. Gronwall's inequalities in [3] show that hazard rates of the form $\exp(\theta^T z)[1 + x\exp(\theta^T z)]^\xi, \xi < -1$, lead to a Volterra equation whose solution is on a finite interval dependent on $\theta$ which may be strictly contained in the support of the withdrawal time distribution. Our results do not apply to this setting.

## 4. Discussion

In Section 2 we discussed properties of the estimate of the unknown transformation under no special regularity on the model representing the "true" distribution $P$ of the data. Examples of Section 3 show that the class of transformation models to which these results apply is quite large and allows hazards of core models to have a variety of shapes.

To estimate the unknown Euclidean $\theta$ parameter, we made the additional assumption that the failure and censoring times are conditionally independent given the covariates and the failure times follow the transformation model. These conditions are sufficient to ensure that the score process is asymptotically unbiased, and the solution to the score equation forms a consistent estimate of the "true" parameter $\theta_0$. However, only first two moment characteristics of certain stochastic integrals are used for this purpose in Section 7, so that the results may also be valid under different assumptions on the true distribution $P$.

We also showed that the class of M-estimators includes a special choice corresponding to an estimate whose asymptotic variance is equal to the inverse of the asymptotic variance of the score function $\sqrt{n}U_n(\theta_0)$. In [17] we show that this estimate is asymptotically efficient. Therein we discuss alternative ad hoc estimators of the unknown transformation and consider a larger class of M-estimators, allowing to adjust common inefficient estimates of the $\theta$ parameter to efficient one-step MLE estimates. Note that asymptotic variance of an M-estimator is usually of a "sandwich" form : As. var $\sqrt{n}(\hat\theta - \theta_0) = A^{-1}(\theta_0)$ As. var$\sqrt{n}U(\theta_0)[A^T(\theta_0)]^{-1}$, where $A(\theta)$ is the limit in probability of the derivative of $U_n(\theta)$ with respect to $\theta$. However, it is quite common that estimators derived from conditional likelihoods of type (2.6) satisfy $A(\theta_0) = $ As. var $\sqrt{n}U_n(\theta_0)$ but are inefficient, so that results of Proposition 2.3 do not imply asymptotic efficiency of the corresponding estimate of the parameter $\theta$.

The proofs of Propositions 2.1 and 2.2 are based on empirical and U-process techniques and are given in Sections 6 and 7. The next section collects some auxiliary results. The proof of consistency and weak convergence of the estimate of the unknown transformation relies also on Gronwall's inequalities collected in Section 9. The proof of Proposition 2.3 uses Fredholm determinant formula for resolvents of linear integral equations [25].

## 5. Some auxiliary results

We denote by $P_n = n^{-1}\sum_{i=1}^n \varepsilon_{X_i,\delta_i,Z_i}$ the empirical measure corresponding to a sequence of $n$ iid observations $(X_i, \delta_i, Z_i)$ representing withdrawal times, censoring indicators and covariates. Set $N.(t) = n^{-1}\sum_{i=1}^n 1(X_i \leq t, \delta_i = 1)$, $Y.(t) =$ $n^{-1}\sum_{i=1}^{n} 1(X_i \geq t)$ and Further, let $\|\cdot\|$ be the supremum norm in the set $\ell^\infty([0,\tau]\times\Theta)$, and let $\|\cdot\|_\infty$ be the supremum norm $\ell^\infty([0,\tau])$. We assume that the point $\tau$ satisfies Condition 2.0(iii).

Define

$$R_n(t,\theta) = \int_0^t \frac{N_.(du)}{S(\Gamma_\theta(u-),\theta,u)} - \int_0^t \frac{EN(du)}{s(\Gamma_\theta(u-),\theta,u)},$$

$$R_{pn}(t,\theta) = \int_0^t h(\Gamma_\theta(u-),\theta,u)[N_.-EN](du), \quad p=5,6,$$

$$R_{pn}(t,\theta) = \int_0^t H_n(\Gamma_\theta(u-),\theta,u)N(du)$$
$$\qquad - \int_0^t h(\Gamma_\theta(u-),\theta,u)EN(du), \quad p=7,8,$$

$$R_{9n}(t,\theta) = \int_{[0,t)} EN(du) |\int_{(u,t]} \mathcal{P}_\theta(u,w) R_{5n}(dw,\theta)|,$$

$$R_{10n}(t,\theta) = \int_0^t \sqrt{n} R_n(u-,\theta) R_{5n}(du,\theta),$$

$$R_{pn}(t,\theta) = \int_0^t H_n(\Gamma_{n\theta}(u-),\theta,u) N_.(du)$$
$$\qquad - \int_0^t h(\Gamma_\theta(u-),\theta,u) EN(du), \quad p=11,12,$$

$$B_{pn}(t,\theta) = \int_0^t F_{pn}(u,\theta) R_n(du,\theta), \quad p=1,2,$$

where $\mathcal{P}_\theta(u,w)$ is given by (2.3). In addition, $H_n = K'_n$ for $p=7$ or $p=11$, $H_n = \dot{K}_n$ for $p=8$ or $p=12$, $h=k'$ for $p=5,7$ or $p=11$, and $h=\dot{k}$ for $p=6,8$ or 12. Here $k' = -[s'/s^2]$, $\dot{k} = -[\dot{s}/s^2]$, $K'_n = -[S'/S^2]$, $\dot{K}_n = -[\dot{S}/S^2]$. Further, set $F_{1n}(u,\theta) = [\dot{S} - \bar{e}S](\Gamma_\theta(u),\theta,u)$ and $F_{2n}(u,\theta) = [S' - eS'](\Gamma_\theta(u),\theta,u)$.

**Lemma 5.1.** *Suppose that Conditions* 2.0 *and* 2.1 *are satisfied.*

(i) $\sqrt{n} R_n(t,\theta)$ *converges weakly in* $\ell^\infty([0,\tau]\times\Theta)$ *to a mean zero Gaussian process $R$ whose covariance function is given below.*

(ii) $\|R_{pn}\| \to 0$ *a.s., for* $p=5,\ldots,12$.

(iii) $\sqrt{n}\|B_{pn}\| \to 0$ *a.s. for* $p=1,2$.

(iv) *The processes* $V_n(\Gamma_\theta(t-),\theta,t)$ *and* $V_n(\Gamma_\theta(t),\theta,t)$, *where* $V_n = S/s - 1$ *satisfy* $\|V_n\| = O(b_n)$ *a.s. In addition,* $\|V_n\| \to 0$ *a.s. for* $V_n = [S'-s']/s, [S''-s'']/s, [\dot{S}-\dot{s}]/s, [\ddot{S}-\ddot{s}]/s$ *and* $[\dot{S}'-\dot{s}']/s$.

*Proof.* The Volterra identity (2.2), which defines $\Gamma_\theta$ as a parameter dependent on $P$, is used in the foregoing to compute the asymptotic covariance function of the process $R_{1n}$. In Section 6 we show that the solution to the identity (2.2) is unique and, for some positive constants $d_0, d_1, d_2$, we have

(5.1)
$$\Gamma_\theta(t) \leq d_0 A_P(t), \quad |\Gamma_\theta(t) - \Gamma_{\theta'}(t)| \leq |\theta-\theta'| d_1 \exp[d_2 A_P(t)],$$
$$|\Gamma_\theta(t) - \Gamma_\theta(t')| \leq d_0 |A_P(t) - A_P(t')|$$
$$\leq \frac{d_0}{E_P Y(\tau)} P(X \in (t\wedge t', t\vee t'], \delta=1),$$

with similar inequalities holding for the left continuous version of $\Gamma_\theta = \Gamma_{\theta,P}$. Here $A_P(t)$ is the cumulative hazard function corresponding to observations $(X,\delta)$.

To show part (i), we use the quadratic expansion, similar to the expansion of the ordinary Aalen–Nelson estimator in [19]. We have $R_n = \sum_{j=1}^{4} R_{jn}$,

$$R_{1n}(t,\theta) = \frac{1}{n} \sum_{i=1}^{n} \int_0^t \left[ \frac{N_i(du)}{s(\Gamma_\theta(u-),\theta,u)} - \frac{S_i}{s^2}(\Gamma_\theta(u-),\theta,u) EN(du) \right]$$

$$= \frac{1}{n} \sum_{i=1}^{n} R_{1n}^{(i)}(t,\theta),$$

$$R_{2n}(t,\theta) = \frac{-1}{n^2} \sum_{i \neq j} \int_0^t \left( \frac{S_i - s}{s^2} \right) (\Gamma_\theta(u-),\theta,u)[N_j - EN_j](du),$$

$$R_{3n}(t,\theta) = \frac{-1}{n^2} \sum_{i=1}^{n} \int_0^t \left( \frac{S_i - s}{s^2} \right) (\Gamma_\theta(u-),\theta,u)[N_i - EN_i](du),$$

$$R_{4n}(t,\theta) = \int_0^t \left( \frac{S - s}{s} \right)^2 (\Gamma_\theta(u-),\theta,u) \frac{N_\cdot(du)}{S(\Gamma_\theta(u-),\theta,u)},$$

where $S_i(\Gamma_\theta(u-),\theta,u) = Y_i(u)\alpha(\Gamma_\theta(u-),\theta,Z_i)$.

The term $R_{3n}$ has expectation of order $O(n^{-1})$. Using Conditions 2.1, it is easy to verify that $R_{2n}$ and $n[R_{3n} - ER_{3n}]$ form canonical U-processes of degree 2 and 1 over Euclidean classes of functions with square integrable envelopes. We have $\|R_{2n}\| = O(b_n^2)$ and $n\|R_{3n} - ER_{3n}\| = O(b_n)$ almost surely, by the law of iterated logarithm for canonical U processes [1]. The term $R_{4n}$ can be bounded by $\|R_{4n}\| \leq \|[S/s] - 1\|^2 m_1^{-1} A_n(\tau)$. But for a point $\tau$ satisfying Condition 2.0(iii), we have $A_n(\tau) = A(\tau) + O(b_n)$ a.s. Therefore part (iv) below implies that $\sqrt{n}\|R_{4n}\| \to 0$ a.s.

The term $R_{1n}$ decomposes into the sum $R_{1n} = R_{1n;1} - R_{1n;2}$, where

$$R_{1n;1}(t,\theta) = \frac{1}{n} \sum_{i=1}^{n} \int_0^t \frac{N_i(du) - Y_i(u)A(du)}{s(\Gamma_\theta(u-),\theta,u)},$$

$$R_{1n;2}(t,\theta) = \int_0^t G(u,\theta) C_\theta(du)$$

and $G(t,\theta) = [S(\Gamma_\theta(u-),\theta,u) - s(\Gamma_\theta(u-),\theta,u)Y_\cdot(u)/EY(u)]$. The Volterra identity (2.2) implies

$$n\text{cov}(R_{1n;1}(t,\theta), R_{1n;1}(t',\theta')) = \int_0^{t \wedge t'} \frac{[1 - A(\Delta u)]\Gamma_\theta(du)}{s(\Gamma_{\theta'}(u-),\theta',u)},$$

$n\text{cov}(R_{1n;1}(t,\theta), R_{1n;2}(t',\theta'))$

$$= \int_0^t \int_0^{u \wedge t'} E[\alpha(\Gamma_{\theta'}(v-), Z, \theta'|X = u, \delta = 1] C_{\theta'}(dv)\Gamma_\theta(du)$$

$$- \int_0^t \int_0^{u \wedge t'} E\alpha(\Gamma_{\theta'}(v-), Z, \theta'|X \geq u]] C_{\theta'}(dv)\Gamma_\theta(du),$$

$n\text{cov}(R_{1n;2}(t,\theta), R_{1n;2}(t',\theta'))$

$$= \int_0^t \int_0^{t' \wedge u} f(u,v,\theta,\theta') C_\theta(du) C_{\theta'}(dv)$$

$$+ \int_0^{t'} \int_0^{t \wedge v} f(v,u,\theta',\theta) C_\theta(du) C_{\theta'}(dv)$$

$$- \int_0^{t \wedge t'} f(u,u,\theta,\theta') C_{\theta'}(\Delta u) C_\theta(du),$$

where $f(u, v, \theta, \theta') = EY(u)\text{cov}(\alpha(\Gamma_\theta(u-), \theta, Z), \alpha(\Gamma_{\theta'}(v-), \theta', Z)|X \geq u)$. Using CLT and Cramer-Wold device, the finite dimensional distributions of $\sqrt{n}R_{1n}(t, \theta)$ converge in distribution to finite dimensional distributions of a Gaussian process. The process $R_{1n}$ can be represented as $R_{1n}(t, \theta) = [P_n - P]h_{t,\theta}$, where $\mathcal{H} = \{h_{t,\theta}(x, d, z) : t \leq \tau, \theta \in \Theta\}$ is a class of functions such that each $h_{t,\theta}$ is a linear combination of 4 functions having a square integrable envelope and such that each is monotone with respect to $t$ and Lipschitz continuous with respect to $\theta$. This is a Euclidean class of functions [29] and $\{\sqrt{n}R_{1n}(t, \theta) : \theta \in \Theta, t \leq \tau\}$ converges weakly in $\ell^\infty([0, \tau] \times \Theta)$ to a tight Gaussian process. The process $\sqrt{n}R_{1n}(t, \theta)$ is asymptotically equicontinuous with respect to the variance semimetric $\rho$. The function $\rho$ is continuous, except for discontinuity hyperplanes corresponding to a finite number of discontinuity points of $EN$. By the law of iterated logarithm [1], we also have $\|R_{1n}\| = O(b_n)$ a.s.

**Remark 5.1.** Under Condition 2.2, we have the identity

$$n\text{cov}(R_{1n;2}(t, \theta_0), R_{1n;2}(t', \theta_0))$$
$$= \sum_{p=1}^{2} n\text{cov}(R_{1n;p}(t, \theta_0, R_{1n;3-p}(t', \theta_0))$$
$$- \int_{[0, t \wedge t']} EY(u)\text{var}(\alpha(\Gamma_{\theta_0}(u-)|X \geq u)C_{\theta_0}(\Delta u)C_{\theta_0}(du).$$

Here $\theta_0$ is the true parameter of the transformation model. Therefore, using the assumption of continuity of the $EN$ function and adding up all terms, $n\text{cov}(R_{1n}(t, \theta_0), R_{1n}(t', \theta_0)) = n\text{cov}(R_{1n;1}(t, \theta_0), R_{1n;1}(t', \theta_0)) = C_{\theta_0}(t \wedge t')$.

Next set $b_\theta(u) = h(\Gamma_\theta(u-), \theta, u)$, $h = k'$ or $h = \dot{h}$. Then $\int_0^t b_\theta(u)N_\cdot(du) = P_n f_{t,\theta}$, where $f_{t,\theta} = 1(X \leq t, \delta = 1)h(\Gamma_\theta(X \wedge \tau-), \theta, X \wedge \tau-)$. The conditions 2.1 and the inequalities (5.1) imply that the class of functions $\{f_{t,\theta} : t \leq \tau, \theta \in \Theta\}$ is Euclidean for a bounded envelope, for it forms a product of a VC-subgraph class and a class of Lipschitz continuous functions with a bounded envelope. The almost sure convergence of the terms $R_{pn}, p = 5, 6$ follows from Glivenko–Cantelli theorem [29].

Next, set $b_\theta(u) = k'(\Gamma_\theta(u-), \theta, u)$ for short. Using Fubini theorem and $|\mathcal{P}_\theta(u, w)| \leq \exp[\int_u^w |b_\theta(s)|EN(ds)]$, we obtain

$$R_{9n}(t, \theta) \leq \int_{(0,t)} EN(du)|R_{5n}(t, \theta) - R_{5n}(u, \theta)|$$
$$+ \int_{(0,t)} EN(du)|\int_{(u,t]} \mathcal{P}_\theta(u, s-)b_\theta(s)EN(ds)[R_{5n}(t, \theta) - R_{5n}(s, \theta)]|$$
$$\leq 2\|R_{5n}\| \int_{[0,t)} EN(du)[1 + \int_{(u,t]} |\mathcal{P}(u, w-)||b_\theta(w)|EN(dw)]$$
$$\leq 2\|R_{5n}\| \int_0^\tau EN(du) \exp[\int_{(u,\tau]} |b_\theta|(s)EN(ds)] \to 0 \quad \text{a.s.}$$

uniformly in $t \leq \tau, \theta \in \Theta$.

Further, we have $R_{10n}(t,\theta) = \sqrt{n} \sum_{p=1}^{4} R_{10n;p}(t,\theta)$, where

$$R_{10n;p}(t;\theta) = \int_0^t R_{pn}(u-,\theta) R_{5n}(du;\theta) =$$
$$= \int_0^t R_{pn}(u-;\theta) k'(\Gamma_\theta(u-),\theta,u)[N_\cdot - EN](du).$$

We have $\|\sqrt{n} R_{10n;p}\| = O(1) \sup_{\theta,t} |\sqrt{n} R_{pn}(u-,\theta)| \to 0$ a.s. for $p = 2, 3, 4$. Moreover, $\sqrt{n} R_{10n;1}(t;\theta) = \sqrt{n} R_{10n;11}(t;\theta) + \sqrt{n} R_{10n;12}(t;\theta)$, where $R_{10n;11}$ is equal to

$$n^{-2} \sum_{i \neq j} \int_0^t R_{1n}^{(i)}(u-,\theta) k'(\Gamma_\theta(u-),\theta,u)[N_j - EN_j](du),$$

while $R_{10n;12}(t,\theta)$ is the same sum taken over indices $i = j$. These are U-processes over Euclidean classes of functions with square integrable envelopes. By the law of iterated logarithm [1], we have $\|R_{10n;11}\| = O(b_n^2)$ and $n\|R_{10n;12} - ER_{10n;12}\| = O(b_n)$ a.s. We also have $ER_{10n;12}(t,\theta) = O(1/n)$ uniformly in $\theta \in \Theta$, and $t \leq \tau$.

The analysis of terms $B_{1n}$ and $B_{2n}$ is quite similar. Suppose that $\ell'(x,\theta) \not\equiv 0$. We have $B_{2n} = \sum_{p=1}^{4} B_{2n;p}$, where in the term $B_{2n;p}$ integration is with respect to $R_{np}$. For $p = 1$, we obtain $B_{2n;1} = B_{2n;11} + B_{2n;12}$, where

$$B_{2n;11}(t,\theta) = \frac{1}{n^2} \int_0^t \sum_{i \neq j} [S'_i - eS_i](\Gamma_\theta(u),\theta,u) R_{1n}^{(j)}(du,\theta),$$

whereas the term $B_{2n;12}$ represents the same sum taken over indices $i = j$. These are U-processes over Euclidean classes of functions with square integrable envelopes. By the law of iterated logarithm [1], we have $\|B_{2n;11}\| = O(b_n^2)$ and $n\|B_{2n;12} - EB_{2n;12}\| = O(b_n)$ a.s. We also have $EB_{2n;12}(t,\theta) = O(1/n)$ uniformly in $\theta \in \Theta$, and $t \leq \tau$. Thus $\sqrt{n}\|B_{2n;1}\| \to 0$ a.s. A similar analysis, leading to U-statistics of degree 1, 2, 3 can be applied to the integrals $\sqrt{n} B_{2n;p}(t,\theta)$, $p = 2, 3$. On the other hand, assumption 2.1 implies that for $p = 4$, we have the bound

$$|B_{2n;4}(t,\theta)| \leq 2 \int_0^\tau \psi(A_2(u-)) \frac{(S-s)^2}{s^2} EN(du)$$
$$\leq O(1) \int_0^\tau \frac{(S-s)^2}{s^2} EN(du),$$

where, under Condition 2.1, the function $\psi$ bounding $\ell'$ is either a constant $c$ or a bounded decreasing function (thus bounded by some c). The right-hand side can further be expanded to verify that $\|\sqrt{n} B_{2n;4}\| \to 0$ a.s. Alternatively, we can use part (iv).

A similar expansion can also be applied to show that $\|R_{7n}\| \to 0$ a.s. Alternatively we have, $|R_{7n}(t,\theta)| \leq \int_0^\tau |K'_n - k'|(\Gamma_\theta(u-),\theta,u) N_\cdot(du) + |R_{5n}(t,\theta)|$ and by part (iv), we have uniform almost sure convergence of the term $R_{7n}$. We also have $|R_{11n} - R_{7n}|(t,\theta) \leq \int_0^\tau O(|\Gamma_{n\theta} - \Gamma_\theta|)(u) N_\cdot(du)$ a.s., so that part (i) implies $\|R_{11n}\| \to 0$ a.s. The terms $R_{8n}$ and $R_{12n}$ can be handled analogously.

Next, $[S/s](\Gamma_\theta(t-),\theta,t) = P_n f_{\theta,t}$, where

$$f_{\theta,t}(x,\delta,z) = 1(x \geq t) \frac{\alpha(\Gamma_\theta(t-),\theta,z)}{EY(u)\alpha(\Gamma_\theta(t-),\theta,Z)} = 1(x \geq t) g_{\theta,t}(z).$$

Suppose that Condition 2.1 is satisfied by a decreasing function $\psi$ and an increasing function $\psi_1$. The inequalities (5.1) and Condition 2.1, imply that $|g_{\theta,t}(Z)| \leq$

$m_2[m_1 E_P Y(\tau)]^{-1}$, $|g_{\theta,t}(Z) - g_{\theta',t}(Z)| \leq |\theta - \theta'| h_1(\tau)$, $|g_{\theta,t}(Z) - g_{\theta,t'}(Z)| \leq [P(X \in [t \wedge t', t \vee t')) + P(X \in (t \wedge t', t \vee t'], \delta = 1)] h_2(\tau)$, where

$$h_1(\tau) = 2m_2[m_1 E_P Y(\tau)]^{-1}[\psi_1(d_0 A_P(\tau)) + \psi(0) d_1 \exp[d_2 A_P(\tau)],$$
$$h_2(\tau) = m_2[m_1 E_P Y(\tau)]^{-2}[m_2 + 2\psi(0)].$$

Setting $h(\tau) = \max[h_1(\tau), h_2(\tau), m_2(m_1 E_P Y(\tau))^{-1}]$, it is easy to verify that the class of functions $\{f_{\theta,t}(x,\delta,z)/h(\tau) : \theta \in \Theta, t \leq \tau\}$ is Euclidean for a bounded envelope. The law of iterated logarithm for empirical processes over Euclidean classes of functions [1] implies therefore that part (iii) is satisfied by the process $V = S/s - 1$. For the remaining choices of the $V$ processes the proof is analogous and follows from the Glivenko–Cantelli theorem for Euclidean classes of functions [29]. □

## 6. Proof of Proposition 2.1

### 6.1. Part (i)

For $P \in \mathcal{P}$, let $A(t) = A_P(t)$ be given by (2.1) and let $\tau_0 = \sup\{t : E_P Y(t) > 0\}$. The condition 2.1 (ii) assumes that there exist constants $m_1 < m_2$ such that the hazard rate $\alpha(x, \theta|z)$ is bounded from below by $m_1$ and from above by $m_2$. Put $A_1 = m_1^{-1} A(t)$ and $A_2(t) = m_2^{-1}(t)$. Then $A_2 \leq A_1$. Further, Condition 2(iii) assumes that the function $\ell(x, \theta, z) = \log \alpha(x, \theta, z)$ has a derivative $\ell'(x, \theta, z)$ with respect to $x$ satisfying $|\ell'(x, \theta, z)| \leq \psi(x)$ for some bounded decreasing function. Suppose that $\psi \leq c$ and define $\rho(t) = \max(c, 1) A_1(t)$. Finally, the derivative $\dot{\ell}(x, \theta, z)$ satisfies $|\dot{\ell}(x, \theta, z)| \leq \psi_1(x)$ for some bounded function or a function that is continuous strictly increasing, bounded at origin and satisfying $\int_0^\infty \psi_1(x)^2 e^{-x} dx < \infty$. Let

$$d = \int_0^\infty \psi_1(x) e^{-x} dx < \infty.$$

In the inequalities (5.1) of Lemma 5.1 we take $d_0 = m_1^{-1}$, $d_1 = \max(1, c)$ and $d_2 = d$.

Let $\mathcal{T} = [0, \tau_0]$ if $\tau_0$ is a discontinuity point of the survival function $E_P Y(t)$, and let $\mathcal{T} = [0, \tau_0)$, if $\tau_0$ is a continuity point of this survival function. Consider the set of functions $\mathcal{X}(P) = \{g : g \text{ monotone increasing}, e^{-g} \in D(\mathcal{T}), g \ll E_P N, A_2 \leq g \leq A_1\}$. Since for each $g \in \mathcal{X}(\mathcal{P})$, the function $e^{-g}$ is a subsurvival function satisfying $\exp[-A_1] \leq \exp[-g] \leq \exp[-A_2]$, we can consider $\mathcal{X}(P)$ as a subset of $D(\mathcal{T})$, endowed with supremum norm. Next, for $\tau < \tau_0$, let $\mathcal{X}(P, \tau) \subset D([0, \tau])$ consist of functions $g \in \mathcal{X}(P)$ restricted to the interval $[0, \tau]$. For fixed $\theta \in \Theta$ and $g \in \mathcal{X}(P, \tau)$, define

$$\Psi_\theta(g)(t) = \int_0^t [EY(u) \alpha(g(u-), \theta, Z)]^{-1} EN(du), \quad 0 \leq t \leq \tau.$$

Using bounds $A_1 \leq g \leq A_2$ it is easy to verify that for fixed $\theta \in \Theta$, $\Psi_\theta$ maps $\mathcal{X}(P, \tau)$ into itself. Since $A_p(0-) = 0$, we have $g(0-) = 0$ and $\Psi_\theta(g)(0-) = 0$ as well.

Consider the equation $\Psi_\theta(g) = g, g(0-) = 0$. Using Helly selection theorem, it is easy to verify that for fixed $\theta \in \Theta$, the operator $\Psi_\theta$ maps $\mathcal{X}(P, \tau)$ into itself, is continuous (with respect to $g$) and has compact range. Since $\mathcal{X}(P, \tau)$ forms a bounded, closed convex set of functions, Schauder's fixed point theorem implies that $\Psi_\theta$ has a fixed point in $\mathcal{X}(P, \tau)$.

To show uniqueness of the solution and its continuity with respect to $\theta$, we consider first the case of continuous $EN(t)$ function. Then $\mathcal{X}(P,\tau) \subset C([0,\tau])$. Define a norm in $C([0,\tau])$ by setting $\|x\|_\rho^\tau = \sup_{t \leq \tau} e^{-\rho(t)}|x(t)|$. Then $\|\cdot\|_\rho^\tau$ is equivalent to the sup norm in $C([0,\tau])$. For $g, g' \in \mathcal{X}(P,\tau)$ and $\theta \in \Theta$, we have

$$\begin{aligned}
|\Psi_\theta(g) - \Psi_\theta(g')|(t) &\leq \int_0^t |g - g'|(u)\psi(A_2(u))A_1(du) \\
&\leq \int_0^t |g - g'|(u)\rho(du) \leq \|g - g'\|_\rho^\tau \int_0^t e^{\rho(u)}\rho(du) \\
&\leq \|g - g'\|_\rho^\tau e^{\rho(t)}(1 - e^{-\rho(\tau)})
\end{aligned}$$

and hence $\|\Psi_\theta(g) - \Psi_{\theta'}(g')\|_\rho^\tau \leq \|g - g'\|_\rho^\tau(1 - e^{-\rho(\tau)})$. For any $g \in \mathcal{X}(P,\tau)$ and $\theta, \theta' \in \Theta$, we also have

$$\begin{aligned}
|\Psi_\theta(g) - \Psi_{\theta'}(g)|(t) &\leq |\theta - \theta'| \int_0^t \psi_1(g(u))A_1(du) \\
&\leq |\theta - \theta'| \int_0^t \psi_1(\rho(u))\rho(du) \\
&\leq |\theta - \theta'|e^{\rho(t)} \int_0^t \psi_1(\rho(u))e^{-\rho(u)}\rho(du) \leq |\theta - \theta'|e^{\rho(t)}d,
\end{aligned}$$

so that $\|\Psi_\theta(g) - \Psi_{\theta'}(g)\|_\rho^\tau \leq |\theta - \theta'|d$. It follows that $\{\Psi_\theta : \theta \in \Theta\}$, restricted to $C[0,\tau])$, forms a family of continuously contracting mappings. Banach fixed point theorem for continuously contracting mappings [24] implies therefore that there exists a unique solution $\Gamma_\theta$ to the equation $\Phi_\theta(g)(t) = g(t)$ for $t \leq \tau$, and this solution is continuous in $\theta$. Since $A(0) = A(0-) = 0$, and the solution is bounded between two multiples of $A(t)$, we also have $\Gamma_\theta(0) = 0$.

Because $\|\cdot\|_\rho^\tau$ is equivalent to the supremum norm in $C[0,\tau]$, we have that for fixed $\tau < \tau_0$, there exists a unique (in sup norm) solution to the equation, and the solution is continuous with respect to $\theta$. It remains to consider the behaviour of these functions at $\tau_0$. Fix $\theta \in \Theta$ again. If $A(\tau_0) < \infty$, then $\Gamma_\theta$ is unique on the whole interval $[0, \tau_0]$ (the preceding argument can be applied to the interval $[0, \tau_0]$). So let us consider the case of $A(\tau) \uparrow \infty$ as $\tau \uparrow \tau_0$. If $\tau^{(1)} < \tau^{(2)} < \tau_0$, then $\mathcal{X}(P, \tau^{(1)}) \subset \mathcal{X}(P, \tau^{(2)})$. Let $\Gamma_\theta^{(p)} \in \mathcal{X}(P, \tau^{(p)}), p = 1, 2$ be the solutions obtained on intervals $[0, \tau^{(1)}]$ and $[0, \tau^{(2)}]$, respectively. Then the function $\Gamma_\theta^{(2)}$ satisfies $\Gamma_\theta^{(2)}(t) = \Gamma_\theta^{(1)}(t)$ for $t \in [0, \tau^{(1)}]$. If $\tau^{(n)} \uparrow \tau_0$, then the inequalities $\exp[-A_1(\tau^{(n)})] \leq \exp[-\Gamma_\theta^{(n)}(\tau^{(n)})] \leq \exp[-A_2(\tau^{(n)})]$ imply $\Gamma_\theta^{(n)}(\tau^{(n)}) \uparrow \infty$. Since this holds for any such sequence $\tau^{(n)}$, there exists a unique locally bounded solution to the equation on the interval $\mathcal{T} = [0, \tau_0)$.

Next let us consider the case of discrete $EN(t)$ with a finite number of discontinuity points. In this case $A_P(\tau_0)$ is bounded and satisfies $A_p(0-) = 0$. Fix $\theta$. Using induction on jumps, it is easy to verify that for any $g \in \mathcal{X}(P, \tau_0)$, we have $\Psi_\theta(g) \in \mathcal{X}(P, \tau_0)$, and for any $g, g' \in \mathcal{X}(P, \tau_0)$, we also have $\Psi_\theta(g) = \Psi_\theta(g')$. Hence $\Psi_\theta^2(g) = \Psi_\theta(g)$. Alternatively, that the solution $\Gamma_\theta$ to the equation $\Psi_\theta(g) = g$, $g(0-) = 0$ is uniquely defined follows also from the recurrent formula $\Gamma_\theta(t) = \Gamma_\theta(t-) + EN(\Delta t)[EY(t)\alpha(\Gamma_\theta(t-), \theta, Z)]^{-1}$, $\Gamma_\theta(0-) = 0$.

For any $g \in \mathcal{X}(P, \tau_0)$ and $\theta, \theta' \in \Theta$, we also have

$$|\Psi_\theta(g) - \Psi_{\theta'}(g)|(t-) \leq |\theta - \theta'| \int_{[0,t)} \psi_1(g(u-))A_1(du)$$

$$\leq |\theta - \theta'| e^{\rho(t)} \int_{[0,t)} \psi_1(\rho(u-)) e^{-\rho(u)} \rho(du)$$

$$\leq |\theta - \theta'| e^{\rho(t-)} d.$$

To see the last inequality, we define

$$\Psi_1(x) = \int_0^x e^{-y} \psi_1(y) dy.$$

Then $\Psi_1(\rho(t)) - \Psi_1(0) = \Sigma \rho(\Delta u) \psi_1(\rho(u^*)) \exp{-\rho(u^*)}$, where the sum extends over discontinuity points less than $t$, and $\rho(u^*)$ is between $\rho(u-)$ and $\rho(u)$. The right-hand side is bounded from below by the corresponding sum $\sum \rho(\Delta u) \psi_1(\rho(u-)) \exp[-\rho(u)]$, because $\psi_1(x)$ is increasing and $\exp(-x)$ is decreasing. Since $\Psi_\theta(g) = \Gamma_\theta$ for any $\theta$, we have $\sup_{t \leq \tau_0} e^{-\rho(t-)}|\Gamma_\theta - \Gamma_{\theta'}|(t-) \leq |\theta - \theta'|d$.

Finally, for both the continuous and discrete case, we have

$$|\Gamma_\theta - \Gamma_{\theta'}|(t) \leq |\Psi_\theta(\Gamma_\theta) - \Psi_\theta(\Gamma_{\theta'})|(t) + |\Psi_\theta(\Gamma_{\theta'}) - \Psi_{\theta'}(\Gamma_{\theta'})|(t) \leq$$

$$\leq \int_0^t |\Gamma_\theta - \Gamma_{\theta'}|(u-)\rho(du) + |\theta - \theta'| \int_0^t \psi_1(\rho(u-))\rho(du),$$

and Gronwall's inequality (Section 9) yields

$$|\Gamma_\theta - \Gamma_{\theta'}|(t) \leq |\theta - \theta'| e^{\rho(t)} \int_{(0,t]} \psi_1(\rho(u-)) e^{-\rho(u-)} \rho(du) \leq d|\theta - \theta'| e^{\rho(t)}.$$

Hence $\sup_{t \leq \tau} e^{-\rho(t)}|\Gamma_\theta - \Gamma_{\theta'}|(t) \leq |\theta - \theta'|d$. In the continuous case this holds for any $\tau < \tau_0$, in the discrete case for any $\tau \leq \tau_0$.

**Remark 6.1.** We have chosen the $\rho$ function as equal to $\rho(t) = \max(c,1)A_1$, where $c$ is a constant bounding the function $\ell'_i(x, \theta)$. Under Condition 2.1, this function may also be bounded by a continuous decreasing function $\psi$. The proof, assuming that $\rho(t) = \int_0^t \psi(A_2(u-))A_1(du)$ is quite similar. In the foregoing we consider the simpler choice, because in Proposition 2.2 we have assumed Condition 2.0(iii). Further, in the discrete case the assumption that the number of discontinuity points is finite is not needed but the derivations are longer.

To show consistency of the estimate $\Gamma_{n\theta}$, we assume now that the point $\tau$ satisfies Condition 2.0(iii). Let $A_n(t)$ be the Aalen–Nelson estimator and set $A_{pn} = m_p^{-1} A_n, p = 1, 2$. We have $A_{2n}(t) \leq \Gamma_{n\theta}(t) \leq A_{1n}(t)$ for all $\theta \in \Theta$ and $t \leq \max(X_i, i = 1, \ldots, n)$. Setting $K_n(\Gamma_{n\theta}(u-), \theta, u) = S(\Gamma_{n\theta}(u-), \theta, u)^{-1}$, we have

$$\Gamma_{n\theta}(t) - \Gamma_\theta(t) = R_n(t, \theta)$$
$$+ \int_{(0,t]} [K_n(\Gamma_{n\theta}(u-), \theta, u) - K_n(\Gamma_\theta(u-), \theta, u)] N_.(du).$$

Hence $|\Gamma_{n\theta}(t) - \Gamma_\theta(t)| \leq |R_n(t,\theta)| + \int_0^t |\Gamma_{n\theta} - \Gamma_\theta|(u-)\rho_n(du)$, where $\rho_n = \max(c, 1)A_{1n}$. Gronwall's inequality implies $\sup_{t, \theta} \exp[-\rho_n(t)]|\Gamma_{n\theta} - \Gamma_\theta|(t) \to 0$ a.s., where the supremum is over $\theta \in \Theta$ and $t \leq \tau$. If $\tau_0$ is a discontinuity point of the survival function $E_P Y(t)$ then this holds for $\tau = \tau_0$.

Next suppose that $\tau_0$ is a continuity point of this survival function, and let $\mathcal{T} = [0, \tau_0)$. We have $\sup_{t \in \mathcal{T}} |\exp[-A_{pn}(t)] - \exp[-A_p(t)]| = o_P(1)$. In addition, for any $\tau < \tau_0$, we have $\exp[-A_{1n}(\tau)] \leq \exp[-\Gamma_{n\theta}(\tau)] \leq \exp[-A_{2n}(\tau)]$. Standard monotonicity arguments imply $\sup_{t \in \mathcal{T}} |\exp[-\Gamma_{n\theta}(t)] - \exp[-\Gamma_\theta(t)]| = o_P(1)$, because $\Gamma_\theta(\tau) \uparrow \infty$ as $\tau \uparrow \infty$.

### 6.2. Part (iii)

The process $\hat{W}(t, \theta) = \sqrt{n}[\Gamma_{n\theta} - \Gamma_\theta](t)$ satisfies

$$\hat{W}(t, \theta) = \sqrt{n} R_n(t, \theta) - \int_{[0,t]} \hat{W}(u-, \theta) N_\cdot(du) b_{n\theta}^*(u),$$

where

$$b_{n\theta}^*(u) = \left[ \int_0^1 (S'/S^2)(\theta, \Gamma_\theta(u-) + \lambda[\Gamma_{n\theta} - \Gamma_\theta](u-), u) d\lambda \right].$$

Define

$$\tilde{W}(t, \theta) = \sqrt{n} R_n(t, \theta) - \int_0^t \tilde{W}(u-, \theta) b_\theta(u) EN(du),$$

where $b_\theta(u) = [s'/s^2](\Gamma_\theta(u), \theta, u)$. We have

$$\tilde{W}(t, \theta) = \sqrt{n} R_n(t, \theta) - \int_0^t \sqrt{n} R_n(u-, \theta) b_\theta(u) EN(du) \mathcal{P}_\theta(u, t)$$

and

$$\hat{W}(t, \theta) - \tilde{W}(t, \theta) = -\int_0^t [\hat{W} - \tilde{W}](u-, \theta) b_{n\theta}^*(u) N_\cdot(du) + \text{rem}(t, \theta),$$

where

$$\text{rem}(t, \theta) = -\int_{[0,t]} \tilde{W}(u-, \theta)[b_{n\theta}^*(u) N_\cdot(du) - b_\theta(u) EN(du)].$$

The remainder term is bounded by

$$\int_0^\tau |\tilde{W}(u-, \theta)||[b_{n\theta}^* - b_\theta](u)| N_\cdot(du) + R_{10n}(t, \theta)$$
$$+ \int_0^{t-} |\sqrt{n} R_n(u-, \theta)||b_\theta(u)| R_{9n}(du, \theta).$$

By noting that $R_{9n}(\cdot, \theta)$ is a nonnegative increasing process, we have $\|\text{rem}\| = o_P(1) + \|R_{10n}\| + O_P(1)\|R_{9n}\| = o_P(1)$. Finally,

$$|\hat{W}(t, \theta) - \tilde{W}(t, \theta)| \leq |\text{rem}(t, \theta)| + \int_0^t |\hat{W} - \tilde{W}|(u-, \theta) \rho_n(du).$$

By Gronwall's inequality (Section 9), we have $\hat{W}(t, \theta) = \tilde{W}(t, \theta) + o_P(1)$ uniformly in $t \leq \tau, \theta \in \Theta$. This verifies that the process $\sqrt{n}[\Gamma_{n\theta} - \Gamma_\theta]$ is asymptotically Gaussian, under the assumption that observations are iid, but Condition 2.2 does not necessarily hold.

### 6.3. Part (ii)

Put

$$\dot{\Gamma}_{n\theta}(t) = \int_0^t \dot{K}_n(\Gamma_{n\theta}(u-), \theta, u) N_\cdot(du) \tag{6.1}$$
$$+ \int_0^t K'_n(\Gamma_{n\theta}(u-), \theta, u) \dot{\Gamma}_{n\theta}(u-) N_\cdot(du),$$

$$\dot{\Gamma}_{\theta}(t) = \int_0^t \dot{k}(\Gamma_{\theta}(u-), \theta, u) EN(du) \tag{6.2}$$
$$+ \int_0^t k'(\Gamma_{\theta}(u-), \theta, u) \dot{\Gamma}_{\theta}(u-) EN(du).$$

Here $\dot{K} = \dot{S}/S^2$, $K' = -S'/S^2$, $\dot{k} = \dot{s}/s^2$ and $k' = -s'/s^2$. Assumption 2.0(iii) implies that $\Gamma_\theta(\tau) \le m_1^{-1}(\tau) < \infty$. For $G = k', \dot{k}$, Conditions 2.1 imply that $\sup_{\theta,t} \int_0^t |G(\Gamma_\theta(u-), \theta, u)| EN(du) < \infty$ so that $\sup_\theta \|\dot{\Gamma}_\theta\|_\infty < \infty$. Uniform consistency of the $\Gamma_{n\theta}$ process implies also that for $G_n = K'_n, \dot{K}_n$, we have

$$\limsup_n \sup_{\theta,t} \int_0^t |G_n(\Gamma_{n\theta}(u-), \theta, u)| N_\cdot(du) < \infty$$

almost surely. Substracting equation (3.2) from (3.3), we get

$$[\dot{\Gamma}_{n\theta} - \dot{\Gamma}_\theta](t) = \Psi_n(\theta, t) + \int_0^t [\dot{\Gamma}_{n\theta} - \dot{\Gamma}_\theta](u-) K'_n(\Gamma_\theta(u-), \theta, u) N_\cdot(du),$$

where

$$\Psi_n(t, \theta) = R_{12n}(t, \theta) + \int_0^t \dot{\Gamma}_\theta(u-) R_{11n}(du, \theta).$$

By Lemma 5.1 and Fubini theorem, we have $\|\Psi_n\| \to 0$ a.s. And using Gronwall's inequality (Section 9), $|\dot{\Gamma}_{n\theta} - \dot{\Gamma}_\theta|(t) \to 0$ a.s. uniformly in t and $\theta$.

Further, consider the remainder term $\mathrm{rem}_n(h, \theta, t) = \Gamma_{n\theta+h}(t) - \Gamma_{n\theta}(t) - h^T \dot{\Gamma}_{n\theta}(t)$ for $\theta, \theta + h \in \Theta$. Set $h_{2n} = \Gamma_{n,\theta+h} - \Gamma_{n,\theta}$. We have

$$\mathrm{rem}_n(h, \theta, t) = h^T \int_0^t \psi_{2n}(h, \theta, u) N_\cdot(du)$$
$$+ \int_0^t \mathrm{rem}_n(h, \theta, u-) \psi_{1n}(h, \theta, u) N_\cdot(du),$$

$$\psi_{1n}(h, \theta, u) = \int_0^1 K'_n(\Gamma_{n\theta}(u-) + \lambda h_{2n}(u-), \theta + \lambda h, u) d\lambda,$$
$$\psi_{2n}(h, \theta, u)$$
$$= \int_0^1 \left[\dot{K}_n(\Gamma_{n\theta}(u-) + \lambda h_{2n}(u-), \theta + \lambda h, u) - \dot{K}_n(\Gamma_{n\theta}(u-), \theta, u)\right] d\lambda$$
$$+ \int_0^1 [K'_n(\Gamma_{n\theta}(u-) + \lambda h_{2n}(u-), \theta + \lambda h, u) - K'_n(\Gamma_{n\theta}(u-), \theta, u)]$$
$$\times \dot{\Gamma}_{n\theta}(u-) d\lambda.$$

We have $\int_0^t |\psi_{1n}(h,\theta,u)| N_{\cdot}(du) \leq \rho_n(t)$ and $\int_0^t |\psi_{2n}(h,\theta,u)| N_{\cdot}(du) \leq h^T \int_0^t B_n(u) \times N_{\cdot}(du)$, for a process $B_n$ with $\limsup_n \int_0^\tau B_n(u) N_{\cdot}(du) = O(1)$ a.s. This follows from condition 2.1 and some elementary algebra. By Gronwall's inequality, $\limsup_n \sup_{t \leq \tau} |\text{rem}_n(h,\theta,t)| = O(|h|^2) = o(|h|)$ a.s. A similar argument shows that if $h_n$ is a nonrandom sequence with $h_n = O(n^{-1/2})$, then $\limsup_n \sup_{t \leq \tau} |\text{rem}_n(h_n, \theta,t)| = O(n^{-1})$ a.s. If $\hat{h}_n$ is a random sequence with $|\hat{h}_n| \xrightarrow{P} 0$, then $\limsup_n \sup_{t \leq \tau} |\text{rem}_n(\hat{h}_n,\theta,t)| = O_p(|\hat{h}_n|^2)$.

### 6.4. Part (iv)

Next suppose that $\theta_0$ is a fixed point in $\Theta$, $EN(t)$ is continuous, and $\hat{\theta}$ is a $\sqrt{n}$-consistent estimate of $\theta_0$. Since $EN(t)$ is a continuous function, $\{\hat{W}(t,\theta): t \leq \tau, \theta \in \Theta\}$ converges weakly to a process $W$ whose paths can be taken to be continuous with respect to the supremum norm. Because $\sqrt{n}[\hat{\theta} - \theta_0]$ is bounded in probability, we have $\sqrt{n}[\Gamma_{n\hat{\theta}} - \Gamma_{\theta_0}] - \sqrt{n}[\hat{\theta} - \theta_0]\dot{\Gamma}_{\theta_0} = \hat{W}(\cdot,\hat{\theta}) + \sqrt{n}[\Gamma_{\hat{\theta}} - \Gamma_{\theta_0} - [\hat{\theta} - \theta_0]\dot{\Gamma}_{\theta_0}] = \hat{W}(\cdot,\hat{\theta}) + O_P(\sqrt{n}|\hat{\theta} - \theta_0|^2) \Rightarrow W(\cdot,\theta_0)$ by weak convergence of the process $\{\hat{W}(t,\theta): t \leq \tau, \theta \in \Theta\}$ and [8].

### 7. Proof of Proposition 2.2

The first part follows from Remark 3.1 and part (iv) of Proposition 2.1. Note that at the true parameter value $\theta = \theta_0$, we have $\sqrt{n}[\Gamma_{n\theta_0} - \Gamma_{\theta_0}](t) = n^{1/2} \int_0^t R_{1n}(du, \theta_0) \times \mathcal{P}_{\theta_0}(u,t) + o_P(1)$, where $R_{1n}$ is defined as in Lemma 5.1,

$$R_{1n}(t,\theta) = \frac{1}{n} \sum_{i=1}^n \int_0^t \frac{M_i(du,\theta)}{s(\Gamma_\theta(u-),\theta,u)}.$$

and $M_i(t,\theta) = N_i(t) - \int_0^t Y_i(u) \alpha(\Gamma_\theta, \theta, Z_i) \Gamma_\theta(du)$.

We shall consider now the score process. Define

$$\tilde{U}_{n1}(\theta) = \frac{1}{n} \sum_{i=1}^n \int_0^\tau \tilde{b}_i(\Gamma_\theta(u), \theta) M_i(dt, \theta),$$

$$\tilde{U}_{n2}(\theta) = \int_0^\tau R_{1n}(du, \theta) \int_{(u,\tau]} \mathcal{P}_\theta(u, v-) r(dv, \theta).$$

Here $\tilde{b}_i(\Gamma_\theta(u), \theta) = \tilde{b}_{i1}(\Gamma_\theta(u), \theta) - \tilde{b}_{i2}(\Gamma_\theta(u), \theta) \varphi_{\theta_0}(t)$ and $\tilde{b}_{1i}(\Gamma_\theta(t), \theta) = \dot{\ell}(\Gamma_\theta(t), \theta, Z_i) - [\dot{s}/s](\Gamma_\theta(t), \theta, t)$, $\tilde{b}_{2i}(\Gamma_\theta(t), \theta) = \ell'(\Gamma_\theta(t), \theta, Z_i) - [s'/s](\Gamma_\theta(t), \theta, t)$. The function $r(\cdot, \theta)$ is the limit in probability of the term $\hat{r}_1(t, \theta)$ given below. Under Condition 2.2, it reduces at $\theta = \theta_0$ to

$$r(\cdot, \theta_0) = -\int_0^t \rho_\varphi(u, \theta_0) EN(du)$$

and $\rho_\varphi(u, \theta_0)$ is the conditional correlation defined in Section 2.3. The terms $\sqrt{n}\tilde{U}_{1n}(\theta_0)$ and $\sqrt{n}\tilde{U}_{2n}(\theta_0)$ are uncorrelated sums of iid mean zero variables and their sum converges weakly to a mean zero normal variable with covariance matrix $\Sigma_{2,\varphi}(\theta_0, \tau)$ given in the statement of Proposition 2.2.

We decompose the process $U_n(\theta)$ as $U_n(\theta) = \hat{U}_n(\theta) + \bar{U}_n(\theta)$, where

$$\hat{U}_n(\theta) = \frac{1}{n}\sum_{i=1}^n \int_0^\tau [b_i(\Gamma_{n\theta}(t),\theta,t) - b_{2i}(\Gamma_{n\theta}(t),\theta,t)\varphi_{\theta_0}(t)]N_i(dt),$$

$$\bar{U}_n(\theta) = -\frac{1}{n}\sum_{i=1}^n \int_0^\tau b_{2i}(\Gamma_{n\theta}(t),\theta,t)[\varphi_{n\theta} - \varphi_{\theta_0}](t)]N_i(dt).$$

We have $\hat{U}_n(\theta) = \sum_{j=1}^3 U_{nj}(\theta)$, where

$$U_{n1}(\theta) = \tilde{U}_{n1}(\theta) + B_{n1}(\tau,\theta) - \int_0^\tau \varphi_{\theta_0}(u) B_{n2}(du,\theta),$$

$$U_{n2}(\theta) = \int_0^\tau [\Gamma_{n\theta} - \Gamma_\theta](t)\hat{r}_1(dt,\theta),$$

$$U_{n3}(\theta) = \int_0^\tau [\Gamma_{n\theta} - \Gamma_\theta](t)\hat{r}_2(dt,\theta).$$

As in Section 2.4, $b_{1i}(x,\theta,t) = \dot{\ell}(x,\theta,Z_i) - [\dot{S}/S](x,\theta,t)$ and $b_{2i}(x,\theta,t) = \ell'(x,\theta,Z_i) - [S'/S](x,\theta,t)$. If $\dot{b}_{pi}$ and $b'_{pi}$ are the derivatives of these functions with respect to $\theta$ and $x$, then

$$\hat{r}_1(s,\theta) = \frac{1}{n}\sum_{i=1}^n \int_0^s [b'_{1i}(\Gamma_\theta(t),\theta,t) - b'_{2i}(\Gamma_\theta(t),\theta,t)\varphi_{\theta_0}(t)]N_i(dt),$$

$$\hat{r}_2(s,\theta) = \frac{1}{n}\sum_{i=1}^n \int_0^s \int_0^1 \hat{r}_{2i}(t,\theta,\lambda)d\lambda N_i(dt),$$

$$\hat{r}_{2i}(t,\theta,\lambda) = [b'_{1i}(\Gamma_\theta(t) + \lambda(\Gamma_{n\theta} - \Gamma_\theta)(t),\theta,t) - b'_{1i}(\Gamma_\theta(t),\theta,t)]$$
$$- [b'_{2i}(\Gamma_\theta(t) + \lambda(\Gamma_{n\theta} - \Gamma_\theta)(t),\theta,t) - b'_{2i}(\Gamma_\theta(t),\theta,t)]\varphi_{\theta_0}(t).$$

We also have $\bar{U}_n(\theta) = U_{n4}(\theta) + U_{n5}(\theta)$, where

$$U_{n4}(\theta) = -\int_0^\tau [\varphi_{n\theta} - \varphi_{\theta_0}](t) B_n(dt,\theta),$$

$$U_{n5}(\theta) = \frac{1}{n}\sum_{i=1}^n \int_0^\tau [\varphi_{n\theta} - \varphi_{\theta_0}](t)[b_{2i}(\Gamma_{n\theta}(u),\theta,u) - b_{2i}(\Gamma_\theta(u),\theta,u)]N_i(dt),$$

$$B_n(t,\theta) = \frac{1}{n}\sum_{i=1}^n \int_0^t b_{2i}(\Gamma_\theta(u),\theta,u)N_i(du)$$
$$= \frac{1}{n}\sum_{i=1}^n \int_0^t \tilde{b}_{2i}(\Gamma_\theta(u),\theta,u)M_i(du,\theta) + B_{2n}(t,\theta).$$

We first show that $\bar{U}_n(\theta_0) = o_P(n^{-1/2})$. By Lemma 5.1, $\sqrt{n}B_{2n}(t,\theta_0)$ converges in probability to 0, uniformly in $t$. At $\theta = \theta_0$, the first term multiplied by $\sqrt{n}$ converges weakly to a mean zero Gaussian martingale. We have $\|\varphi_{n\theta_0} - \varphi_{\theta_0}\|_\infty = o_P(1)$, $\|\varphi_{\theta_0}\|_v < \infty$ and $\limsup_n \|\varphi_{n\theta_0}\|_v < \infty$. Integration by parts, Skorohod–Dudley construction and arguments similar to Lemma A.3 in [7], show that $\sqrt{n}U_{n4}(\theta_0) = o_P(1)$. We also have $\sqrt{n}U_{n5}(\theta_0) = \int_0^\tau O_P(\sqrt{n}|\Gamma_{n\theta_0} - \Gamma_{\theta_0}|(t)|\varphi_{n\theta_0} - \varphi_0|(t)N_\cdot(dt) = o_P(1)$.

We consider now the term $\hat{U}_n(\theta_0)$. We have $\sqrt{n}U_{n3}(\theta_0) = \sqrt{n}\int_0^\tau O_P(|\Gamma_{n\theta_0} - \Gamma_{\theta_0}|(t)^2)N_\cdot(dt) = o_P(1)$. We also have $\|r(\cdot,\theta_0)\|_v < \infty$, $\limsup_n \|\hat{r}(\cdot,\theta_0)\|_v < \infty$ and

$\|[\hat{r}_1 - r](\cdot, \theta_0)\|_\infty = o_P(1)$, so that the same integration by parts argument implies that $\sqrt{n}U_{n2}(\theta_0) = \sqrt{n}\tilde{U}_{n2}(\theta_0) + o_P(1)$. Finally, $\sqrt{n}U_{n1}(\theta_0) = \sqrt{n}\tilde{U}_{n1}(\theta_0) + o_P(1)$, by Lemma 5.1 and Fubini theorem.

Suppose now that $\theta$ varies over a ball $B(\theta_0, \varepsilon_n)$ centered at $\theta_0$ and having radius $\varepsilon_n, \varepsilon_n \downarrow 0, \sqrt{n}\varepsilon_n \to \infty$. It is easy to verify that for $\theta, \theta' \in B(\theta_0, \varepsilon_n)$ we have $U_n(\theta') - U_n(\theta) = -(\theta' - \theta)^T \Sigma_{1n}(\theta_0) + (\theta' - \theta)^T R_n(\theta, \theta')$, where $R_n(\theta, \theta')$ is a remainder term satisfying $\sup\{|R_n(\theta, \theta')| : \theta, \theta' \in B(\theta_0, \varepsilon_n)\} = o_P(1)$. The matrix $\Sigma_{1n}(\theta)$ is equal to the sum $\Sigma_{1n}(\theta) = \Sigma_{11n}(\theta) + \Sigma_{12n}(\theta)$,

$$\Sigma_{11n}(\theta) = \frac{1}{n}\sum_{i=1}^n \int_0^\tau [g_{1i}g_{2i}^T](\Gamma_{n\theta}(u), \theta, u)^T N_i(du),$$

$$\Sigma_{12n}(\theta) = -\frac{1}{n}\sum_{i=1}^n \int_0^\tau [f_i - S_f/S](\Gamma_{n\theta}(u), \theta, u) N_i(du),$$

where $S_f(\Gamma_{n\theta}(u), \theta, u) = n^{-1}\sum_{i=1}^n Y_i(u)[\alpha_i f_i](\Gamma_{n\theta}(u), \theta, u)$ and

$$g_{1i}(\theta, \Gamma_{n\theta}(u), u) = b_{1i}(\Gamma_{n\theta}(u), \theta) - b_{2i}(\Gamma_{n\theta}(u), \theta)\varphi_{\theta_0}(u),$$
$$g_{2i}(\theta, \Gamma_{n\theta}(u), u) = b_{1i}(\Gamma_{n\theta}(u), \theta) + b_{2i}(\Gamma_{n\theta}(u), \theta)\dot{\Gamma}_{n\theta}(u)$$
$$f_i(\theta, \Gamma_{n\theta}(u), u) = \frac{\ddot{\alpha}}{\alpha}(\Gamma_{n\theta}(u), \theta, Z_i) - \frac{\dot{\alpha}'}{\alpha}(\Gamma_{n\theta}(u), \theta, Z_i)\varphi_{\theta_0}(u)^T$$
$$+ \dot{\Gamma}_{n\theta}(u)[\frac{\dot{\alpha}'}{\alpha}(\Gamma_{n\theta}(u), \theta, Z_i)]^T$$
$$+ \frac{\alpha''}{\alpha}(\Gamma_{n\theta}(u), \theta, Z_i)\dot{\Gamma}_{n\theta}(u)\varphi_{\theta_0}(u)^T.$$

These matrices satisfy $\Sigma_{11n}(\theta_0) \to_P \Sigma_{1,\varphi}(\theta_0, \tau)$ and $\Sigma_{12n}(\theta_0) \to_P 0$, and $\Sigma_{1,\varphi}(\theta_0, \tau) = \Sigma_1(\theta_0)$ is defined in the statement of Proposition 2.2. By assumption this matrix is non-singular. Finally, set $h_n(\theta) = \theta + \Sigma_1(\theta_0)^{-1}U_n(\theta)$. It is easy to verify that this mapping forms a contraction on the set $\{\theta : |\theta - \theta_0| \leq A_n/(1-a_n)\}$, where $A_n = |\Sigma_1(\theta_0)^{-1}U_n(\theta_0)| = O_P(n^{-1/2})$ and $a_n = \sup\{|I - \Sigma_1(\theta_0)^{-1}\Sigma_{1n}(\theta_0) + \Sigma_1(\theta_0)^{-1}R_n(\theta, \theta')| : \theta, \theta' \in B(\theta_0, \varepsilon_n)\} = o_P(1)$. The argument is similar to Bickel et al. ([6], p.518), though note that we cannot apply their mean value theorem arguments.

Next consider Condition 2.3(v.2). In this case we have $\hat{U}_n(\theta') - \hat{U}_n(\theta) = -(\theta' - \theta)^T\Sigma_{1n}(\theta_0) + (\theta' - \theta)^T\hat{R}_n(\theta, \theta')$, where $\sup\{|\hat{R}_n(\theta, \theta') : \theta, \theta' \in B(\theta_0, \varepsilon_n)\} = o_P(1)$. In addition, for $\theta \in B_n(\theta_0, \varepsilon)$, we have the expansion $\bar{U}_n(\theta) = [\bar{U}_n(\theta) - \bar{U}_n(\theta_0)] + \bar{U}_n(\theta_0) = o_P(|\theta - \theta_0| + n^{-1/2})$. The same argument as above shows that the equation $\hat{U}_n(\theta)$ has, with probability tending to 1, a unique root in the ball $B(\theta_0, \varepsilon_n)$. But then, we also have $U_n(\hat{\theta}_n) = \hat{U}_n(\hat{\theta}_n) + \bar{U}_n(\hat{\theta}_n) = o_P(|\hat{\theta}_n - \theta_0| + n^{-1/2}) = o_P(O_P(n^{-1/2}) + n^{-1/2}) = o_p(n^{-1/2})$.

Part (iv) can be verified analogously, i.e. it amounts to showing that if $\sqrt{n}[\hat{\theta} - \theta_0]$ is bounded in probability, then the remainder term $\hat{R}_n(\hat{\theta}, \theta_0)$ is of order $o_P(|\hat{\theta} - \theta_0|)$, and $\bar{U}_n(\hat{\theta}) = o_P(|\hat{\theta} - \theta_0| + n^{-1/2})$.

## 8. Proof of Proposition 2.3

Part (i) is verified at the end of the proof. To show part (ii), define

$$D(\lambda) = \sum_{m\geq 0} \frac{(-1)^m}{m!} \lambda^m d_m,$$

$$D(t,u,\lambda) = \sum_{m\geq 0} \frac{(-1)^m}{m!} \lambda^m D_m(t,u).$$

The numbers $d_m$ and the functions $D_m(t,u)$ are given by $d_m = 1$, $D_m(t,u) = k(t,u)$ for $m = 0$. For $m \geq 1$ set

$$d_m = \int \cdots \int_{(s_1,\ldots,s_m) \in (0,\tau]} \det \bar{d}_m(\mathbf{s}) b(ds_1) \cdot \ldots \cdot b(ds_m),$$

$$D_m(t,u) = \int \cdots \int_{(s_1,\ldots,s_m) \in (0,\tau]} \det \bar{D}_m(t,u;\mathbf{s}) b(ds_1) \cdot \ldots \cdot b(ds_m),$$

where for any $\mathbf{s} = (s_1,\ldots,s_m)$, $\bar{d}_m(\mathbf{s})$ is an $m \times m$ matrix with entries $\bar{d}_m(\mathbf{s}) = [k(s_i,s_j)]$, and $\bar{D}_m(t,u;\mathbf{s})$ is an $(m+1) \times (m+1)$ matrix

$$\bar{D}_m(t,u;\mathbf{s}) = \begin{pmatrix} k(t,u), & U_m(t;\mathbf{s}) \\ V_m(\mathbf{s};u), & \bar{d}_m(\mathbf{s}) \end{pmatrix},$$

where $U_m(t;\mathbf{s}) = [k(t,s_1),\ldots,k(t,s_m)]$, $V_m(\mathbf{s};u) = [k(s_1,u),\ldots,k(s_m,u)]^T$. By Fredholm determinant formula [25], the resolvent of the kernel $k$ is given by $\tilde{\Delta}(t,u,\lambda) = D(t,u,\lambda)/D(\lambda)$, for all $\lambda$ such that $D(\lambda) \neq 0$, so that

$$d_m = \int \cdots \int_{\substack{s_1,\ldots,s_m \in (0,\tau] \\ \text{distinct}}} \det \bar{d}_m(\mathbf{s}) b(ds_1) \cdot \ldots \cdot b(ds_m),$$

because the determinant is zero whenever two or more points $s_i, i = 1,\ldots,m$ are equal. By Fubini theorem, the right-hand side of the above expression is equal to

$$\sum_\pi \int \cdots \int_{0<s_1<s_2<\cdots<s_m \leq \tau} \det \bar{d}_m(s_{\pi(1)},\ldots,s_{\pi(m)}) b(ds_1) \cdot \ldots \cdot b(ds_m)$$

$$= m! \int \cdots \int_{0<s_1<s_2<\cdots<s_m \leq \tau} \det \bar{d}_m(s_1,\ldots,s_m) b(ds_1) \cdot \ldots \cdot b(ds_m).$$

The first sum extends over the $m!$ possible permutations $\pi = (\pi(1),\ldots,\pi(m))$ of the index set $\{1,\ldots,m\}$. The second line follows by noting that

$$\det \bar{d}_m(s_{\pi(1)},\ldots,s_{\pi(m)}) = \det \bar{d}_m(s_1,\ldots,s_m)$$

for any such permutation, because the matrix $\bar{d}_m(s_{\pi(1)},\ldots,s_{\pi(m)})$ is symmetric and to rearrange it into the matrix $\bar{d}_m(s_1,\ldots,s_m)$, we need the same number of transpositions of rows and columns. Since the total number of such transpositions is even, the determinants have the same sign. In the same way, the function $D_m(t,u)$ is equal to

$$m! \int \cdots \int_{0<s_1<s_2<\cdots<s_m \leq \tau} \det \bar{D}_m(t,u;s_1,\ldots,s_m) b(ds_1) \cdot \ldots \cdot b(ds_m),$$

so that in both cases it is enough to consider the determinants for ordered sequences $\mathbf{s} = (s_1, \ldots, s_m), s_1 < s_2 < \ldots < s_m$ of points in $(0, \tau]^m$.

For any such sequence $\mathbf{s}$, the matrix $\bar{d}_m(\mathbf{s})$ has a simple pattern:

$$\bar{d}_m(\mathbf{s}) = \begin{pmatrix} c(s_1) & c(s_1) & c(s_1) & \ldots & c(s_1) \\ c(s_1) & c(s_2) & c(s_2) & \ldots & c(s_2) \\ c(s_1) & c(s_2) & c(s_3) & \ldots & c(s_3) \\ \vdots & & & & \vdots \\ c(s_1) & c(s_2) & c(s_3) & \ldots & c(s_m) \end{pmatrix}.$$

We have $\bar{d}_m(\mathbf{s}) = A_m^T C_m(\mathbf{s}) A_m$ where $C_m(\mathbf{s})$ is a diagonal matrix of increments

$$C_m(\mathbf{s}) = \text{diag}\,[c(s_1) - c(s_0), c(s_2) - c(s_1), \ldots c(s_m) - c(s_{m-1})],$$

$(c(s_0) = 0, s_0 = 0)$ and $A_m$ is an upper triangular matrix

$$A_m = \begin{pmatrix} 1 & 1 & \ldots & 1 & 1 \\ 0 & 1 & \ldots & 1 & 1 \\ \vdots & & & & \vdots \\ 0 & 0 & \ldots & 1 & 1 \\ 0 & 0 & \ldots & 0 & 1 \end{pmatrix}.$$

To see this it is enough to note that Brownian motion forms a process with independent increments, and the kernel $k(s, t) = c(s \wedge t)$ is the covariance function of a time transformed Brownian motion.

Apparently, $\det A_m = 1$. Therefore

$$\det \bar{d}_m(\mathbf{s}) = \prod_{j=1}^{m} [c(s_j) - c(s_{j-1})]$$

and

$$\begin{aligned}\det \bar{D}_m(t, u; \mathbf{s}) &= \det \bar{d}_m(\mathbf{s})[c(t \wedge u) - U_m(t; \mathbf{s})[\bar{d}_m(\mathbf{s})]^{-1} V_m(\mathbf{s}; u)] \\ &= \det \bar{d}_m(\mathbf{s})[c(t \wedge u) - U_m(t; \mathbf{s}) A_m^{-1} C_m^{-1}(\mathbf{s})(A_m^T)^{-1} V_m(\mathbf{s}; u)].\end{aligned}$$

The inverse $A_m^{-1}$ is given by Jordan matrix

$$A_m^{-1} = \begin{pmatrix} 1 & -1 & 0 & \ldots & 0 & 0 \\ 0 & 1 & -1 & \ldots & 0 & 0 \\ & \vdots & & & \vdots & \\ 0 & 0 & 0 & \ldots & 1 & -1 \\ 0 & 0 & 0 & \ldots & 0 & 1 \end{pmatrix}$$

and a straightforward multiplication yields

$$\begin{aligned}\det \bar{D}_m(t, u; \mathbf{s}) &= c(t \wedge u) \prod_{j=1}^{m}[c(s_j) - c(s_{j-1})] \\ &\quad - \sum_{i=1}^{m}[c(t \wedge s_i) - c(t \wedge s_{i-1})][c(u \wedge s_i) - c(u \wedge s_{i-1})] \\ &\quad \times \prod_{j=1, j \neq i}^{m}[c(s_j) - c(s_{j-1})].\end{aligned}$$

By noting that the $i$-th summand is zero whenever $t \wedge u < s_{i-1}$ and using induction on $m$, it is easy to verify that for $t \leq u$ the determinant reduces to the sum

$$\det \bar{D}_m(t,u;\mathbf{s}) = 1(t \leq u < s_1)c(t)[c(s_1) - c(u)] \prod_{j=2}^{m} [c(s_j) - c(s_{j-1})]$$

$$+ 1(s_m < t \leq u) \prod_{j=1}^{m} [c(s_j) - c(s_{j-1})][c(t) - c(s_m)]$$

$$+ \sum_{i=1}^{m-1} 1(s_i < t \leq u < s_{i+1}) \left( \prod_{j=1}^{i} [c(s_j) - c(s_{j-1})][c(t) - c(s_i)] \right)$$

$$\times \left( [c(s_{i+1}) - c(u)] \prod_{j=i+2}^{m} [c(s_j) - c(s_{j-1})] \right),$$

where in the last sum, product over an empty set of indices is interpreted as equal to 1. Thus we have a simple expression for the two determinants. Integration with respect to the product measure $b(ds_1) \ldots b(ds_m)$ and induction on $m$ yields also

$$\frac{1}{m!} d_m = \Psi_{0m}(0, \tau),$$

$$\frac{1}{m!} D_m(t, u) = \sum_{l=0}^{m} \Psi_{1l}(0, t \wedge u) \Psi_{0,m-l}(t \vee u, \tau),$$

for $m \geq 0$. The numerator and denominator of the Fredholm determinant formula are bounded functions for any point $\tau$ satisfying the condition 2.0 (iii). For $\lambda = -1$, the ratio $\tilde{\Delta}(t, u, \lambda) = D(t, u, \lambda)/D(\lambda)$ reduces to the function $\tilde{\Delta}$ given in the statement of part (ii). Using monotonicity of the functions $\Psi_0$ with respect to the length of the interval $(s, t]$, we also have $\tilde{\Delta}(t, u) \leq \Psi_0(0, \tau) c(t \wedge u)$. If $\tau_0$ is a continuity point of $EY(t)$ and $\kappa(\tau_0) < \infty$, then the denominator is bounded, and the inequality is satisfied for any $u, t < \tau_0$.

Parts (iii) and (iv) are easy to verify using this last observation. For example, if $\kappa(\tau_0) < \infty$ then for any $\tilde{\eta} \in L_2(b)$, the Fredholm equation has a unique solution $\tilde{\psi}$ and

$$\|\tilde{\psi}\|_2 \leq \|\tilde{\eta}\|_2 + \left[ \int_0^{\tau_0} \left[ \int_0^{\tau_0} \tilde{\Delta}(t, u) b(du) \tilde{\eta}(u) \right]^2 b(dt) \right]^{1/2}.$$

By Cauchy–Schwarz inequality and monotone convergence, the second term is bounded by

$$\|\tilde{\eta}\|_2 \left[ \int_0^{\tau_0} \int_0^{\tau_0} \tilde{\Delta}^2(t, u) b(du) b(dt) \right]^{1/2}$$

$$\leq \|\tilde{\eta}\|_2 \Psi_0(0, \tau_0) \left[ \int_0^{\tau_0} \int_0^{\tau_0} c(t \wedge u) b(du) b(dt) \right]^{1/2}$$

$$\leq \|\tilde{\eta}\|_2 \Psi_0(0, \tau_0) \left[ \int_0^{\tau_0} \int_0^{\tau_0} c(t) c(u) b(du) b(dt) \right]^{1/2} = \|\tilde{\eta}\|_2 \Psi_0(0, \tau_0) \kappa(\tau_0).$$

Part (v) follows from part (iv) and Lemma 2.1. Part (vi) can be verified using straightforward but laborious algebra.

Part (i). For $u > s$, set $c((s, u]) = c(u) - c(s)$. The $n$-the term of the series $\Psi_{0n}(s,t)$ is given by the multiple integral

$$\int_{s<s_1<\cdots<s_n\leq t} c((s,s_1])b(ds_1)c((s_1,s_2])b(ds_2)\cdots c((s_{n-1},s_n])b(ds_n),$$

and satisfies $\Psi_{0n}(s,t) \leq (n!)^{-1}[I(s,t)]^n$, $I(s,t) = \int_s^t c((s,u])b(u)$. The integral $I(s,t)$ is increasing with the width of the interval $(s,t]$ and is bounded by $\kappa(\tau)$. Thus $\Psi_0(s,t) \leq \exp\kappa(\tau) < \infty$ for all $0 < s < t \leq \tau$. If in addition $\kappa(\tau_0) < \infty$, then $\Psi_0(s,t)$ is bounded for all $0 < s < t \leq \tau_0$. In both circumstances, this implies that the remaining interval functions $\Psi_j(s,t)$, $0 < s \leq t \leq \tau$ are finite for any point $\tau$ satisfying Condition 2.0(iii), and monotonically increasing with the size of the interval.

While the identities can be verified by applying Fubini to each term of the $\Psi_j, j = 0, 1, 2, 3$ series, the following provides an interpretation in terms of linear Volterra equations. First, it is easy to see that the "odd" functions satisfy

$$\Psi_1(s,t) = c((s,t]) + \int_{(s,t]} c((s,u])b(du)\Psi_1(u,t)$$
$$= c((s,t]) + \int_{(s,t]} \Psi_1(s,u)b(du)c((u,t]),$$
$$\Psi_3(s,t) = b([s,t)) + \int_{[s,t)} b([s,u))c(du)\Psi_3(u,t)$$
$$= b([s,t)) + \int_{[s,t)} \Psi_3(s,u)c(du)b([u,t)),$$

so that they form resolvents of linear Volterra equations. The "even" functions $\Psi_0$ and $\Psi_2$ satisfy such equations

$$\Psi_0(s,t) = 1 + \int_{(s,t]} c((s,u])b(du)\Psi_0(u,t)$$
$$= 1 + \int_{(s,t]} \Psi_0(s,u-)c(du)b([u,t]),$$
$$\Psi_2(s,t) = 1 + \int_{[s,t)} b([s,u))c(du)\Psi_2(u,t)$$
$$= 1 + \int_{[s,t)} \Psi_2(s,u+)b(du)c((u,t)).$$

With fixed $t$, the equations

$$h_1(s,t) - \int_{(s,t]} c((s,u])b(du)h_1(u,t) = g_1(s,t),$$
$$h_3(s,t) - \int_{(s,t]} b([s,u))c(du)h_3(u,t) = g_3(s,t),$$

have unique solutions

$$h_1(s,t) = g_1(s,t) + \int_{(s,t]} \Psi_1(s,u)b(du)g_1(u,t),$$
$$h_3(s,t) = g_3(s,t) + \int_{[s,t)} \Psi_3(s,u)c(du)g_3(u,t).$$

The first pair of equations for $\Psi_0$ and $\Psi_2$ in part (i) follows by setting $g_1(s,t) = 1 = g_3(s,t)$. With $s$ fixed, the equations

$$\bar{h}_1(s,t) - \int_{[s,t)} \bar{h}_1(s,u+)b(du)c_1((u,t)) = \bar{g}_1(s,t),$$

$$\bar{h}_3(s,t) - \int_{(s,t]} \bar{h}_3(s,u-)c(du)b_3([u,t]) = \bar{g}_3(s,t),$$

have solutions

$$\bar{h}_1(s,t) = \bar{g}_1(s,t) + \int_{[s,t)} \bar{g}_1(s,u+)b(du)\Psi_1(u,t-),$$

$$\bar{h}_3(s,t) = \bar{g}_3(s,t) + \int_{(s,t]} \bar{g}_3(s,u-)c(du)\Psi_3(u,t+).$$

The second pair of equations for $\Psi_0$ and $\Psi_2$ in part (i) follows by setting $\bar{g}_1(s,t) \equiv 1 \equiv \bar{g}_3(s,t)$. Next, the "odd" functions can be represented in terms of "even" functions using Fubini.

## 9. Gronwall's inequalities

Following Gill and Johansen [18], recall that if $b$ is a cadlag function of bounded variation, $\|b\|_v \leq r_1$ then the associated product integral $\mathcal{P}(s,t) = \prod_{(s,t]}(1+b(du))$ satisfies the bound $|\mathcal{P}(s,t)| \leq \prod_{(s,t]}(1+\|b\|_v(dw)) \leq \exp\|b\|_v(s,t]$ uniformly in $0 < s < t \leq \tau$. Moreover, the functions $s \to \mathcal{P}(s,t), s \leq t \leq \tau$ and $t \to \mathcal{P}(s,t), t \in (s,\tau]$ are of bounded variation with variation norm bounded by $r_1 e^{r_1}$.

The proofs use the following consequence of Gronwall's inequalities in Beesack [3] and Gill and Johansen [18]. If $b$ is a nonnegative measure and $y \in D([0,\tau])$ is a nonnegative function then for any $x \in D([0,\tau])$ satisfying

$$0 \leq x(t) \leq y(t) + \int_{(0,t]} x(u-)b(du), \quad t \in [0,\tau],$$

we have

$$0 \leq x(t) \leq y(t) + \int_{(0,t]} y(u-)b(du)\mathcal{P}(u,t), \quad t \in [0,\tau].$$

Pointwise in $t$, $|x(t)|$ is bounded by

$$\max\{\|y\|_\infty, \|y^-\|_\infty\}[1 + \int_{(0,t]} b(du)\mathcal{P}(u,t)] \leq \{\|y\|_\infty, \|y^-\|_\infty\}\exp[\int_0^t b(du)].$$

We also have $\|e^{-b}|x|\|_\infty \leq \max\{\|y\|_\infty, \|y^-\|_\infty\}$. Further, if $0 \not\equiv y \in D([0,\tau])$ and $b$ is a function of bounded variation then the solution to the linear Volterra equation

$$x(t) = y(t) + \int_0^t x(u-)b(du)$$

is unique and given by

$$x(t) = y(t) + \int_{(0,t]} y(u-)b(du)\mathcal{P}(u,t).$$

168    D. M. Dabrowska

We have $|x(t)| \leq \max\{\|y\|_\infty, \|y^-\|_\infty\} \exp \int_0^t d\|b\|_v$ and $\|\exp[-\int d\|b\|_v]|x|\|_\infty \leq \max\{\|y\|_\infty, \|y^-\|_\infty\}$. If $y_\theta(t)$, and $b_\theta(t) = \int_0^t k_\theta(u) n(du)$ are functions dependent on a Euclidean parameter $\theta \in \Theta \subset R^d$, and $|k_\theta|(t) \leq k(t)$, then these bounds hold pointwise in $\theta$ and

$$\sup_{\substack{t \leq \tau \\ \theta \in \Theta}} \{\exp[-\int_0^t k(u)n(du)]|x_\theta(t)|\} \leq \max\{\sup_{\substack{u \leq \tau \\ \theta \in \Theta}} |y_\theta|(u), \sup_{\substack{u \leq \tau \\ \theta \in \Theta}} |y_\theta(u-)|\}.$$

## Acknowledgement

The paper was presented at the First Erich Leh–mann Symposium, Guanajuato, May 2002. I thank Victor Perez Abreu and Javier Rojo for motivating me to write it. I also thank Kjell Doksum, Misha Nikulin and Chris Klaassen for some discussions. The paper benefited also from comments of an anonymous reviewer and the Editor Javier Rojo.